\documentclass[11pt]{amsart}

\usepackage{amsmath,amsfonts,amsthm,amsopn,amssymb}
\usepackage{verbatim,wasysym,mathrsfs}
\usepackage[left=2.6cm,right=2.6cm,top=3.1cm,bottom=3.1cm]{geometry}
\usepackage{cite,marginnote}
\usepackage{color,enumitem,graphicx}
\usepackage[colorlinks=true,urlcolor=blue, citecolor=red,linkcolor=blue,
linktocpage,pdfpagelabels, bookmarksnumbered,bookmarksopen]{hyperref}
\usepackage[english]{babel}
\usepackage[T1]{fontenc}
\usepackage{lmodern}
\usepackage[hyperpageref]{backref}

\numberwithin{equation}{section}

\makeindex
\newtheorem{theorem}{Theorem}[section]
\newtheorem{proposition}[theorem]{Proposition}
\newtheorem{lemma}[theorem]{Lemma}
\newtheorem{remark}[theorem]{Remark}
\newtheorem{example}[theorem]{Example}
\newtheorem{corollary}[theorem]{Corollary}
\newtheorem{definition}[theorem]{Definition}

\newcommand{\s}{\section}

\newcommand{\R}{\mathbb R}
\newcommand{\al}{\alpha}

\newcommand{\e}{\varepsilon}

\newcommand{\bt}{\begin{theorem}}
\newcommand{\et}{\end{theorem}}
\newcommand{\bl}{\begin{lemma}}
\newcommand{\el}{\end{lemma}}
\newcommand{\bd}{\begin{definition}}
\newcommand{\ed}{\end{definition}}
\newcommand{\bc}{\begin{corollary}}
\newcommand{\ec}{\end{corollary}}
\newcommand{\bp}{\begin{proof}}
\newcommand{\ep}{\end{proof}}
\newcommand{\bx}{\begin{example}}
\newcommand{\ex}{\end{example}}
\newcommand{\bi}{\begin{exercise}}
\newcommand{\ei}{\end{exercise}}
\newcommand{\bo}{\begin{proposition}}
\newcommand{\eo}{\end{proposition}}
\newcommand{\br}{\begin{remark}}
\newcommand{\er}{\end{remark}}
\newcommand{\be}{\begin{equation}}
\newcommand{\ee}{\end{equation}}
\newcommand{\ba}{\begin{align}}
\newcommand{\ea}{\end{align}}
\newcommand{\bn}{\begin{enumerate}}
\newcommand{\en}{\end{enumerate}}
\newcommand{\bg}{\begin{align*}}
\newcommand{\bcs}{\begin{cases}}
\newcommand{\ecs}{\end{cases}}

\def\R{\mathbb R}

\def\Proof{\noindent{\bf Proof}\quad}
\def\qed{\hfill$\square$\smallskip}
\makeatletter
\def\@makefnmark{}
\makeatother

\newcommand{\bean}{\begin{eqnarray*}}
\newcommand{\eean}{\end{eqnarray*}}

\renewcommand{\triangle}{\Delta}
\renewcommand{\epsilon}{\varepsilon}

\title[  A new type of nodal solutions]
{ A new type of nodal solutions to singularly perturbed elliptic equations with supercritical growth}

\author[Z. S. Liu]{Zhisu Liu}
\author[J. C. Wei]{Juncheng Wei}
\author[J. J. Zhang]{Jianjun Zhang}

\address[Z. S. Liu]{\newline\indent School of Mathematics and Physics,
\newline\indent
China University of Geosciences,
\newline\indent
Wuhan, Hubei, 430074, PR China}
\email{\href{mailto:liuzhisu@cug.edu.cn}{liuzhisu@cug.edu.cn}}

\address[J. C. Wei]{\newline\indent
Department of Mathematics,
\newline\indent
University of British Columbia,
\newline\indent
Vancouver V6T 1Z2, Canada}
\email{\href{mailto:jcwei@math.ubc.ca}{jcwei@math.ubc.ca}}

\address[J. J. Zhang]{\newline\indent
College of Mathematics and Statistics,
\newline\indent
Chongqing Jiaotong University,
\newline\indent
 Chongqing 400074, PR China }
\email{\href{mailto:zhangjianjun09@tsinghua.org.cn}{zhangjianjun09@tsinghua.org.cn}}

\thanks{(1) Corresponding author: \texttt{jcwei@math.ubc.ca}}

\thanks{(2) The research of J. Wei is partially supported by NSERC of Canada. Z. Liu was supported by the NSFC (No.~11701267), the Hunan Natural Science Excellent Youth Fund (No.~2020JJ3029) and the Fundamental Research Funds for the Central Universities, China University of Geosciences (Wuhan, No. CUG2106211; CUGST2). J. J. Zhang is supported by NSFC(No.11871123).}

\date{\today}

\begin{document}

\begin{abstract}
In this paper, we aim to investigate the following class of singularly perturbed elliptic problems
$$
\left\{
  \begin{array}{ll}
  \displaystyle  -\varepsilon^2\triangle {u}+|x|^\eta u =|x|^\eta f(u)&  \mbox{in}\,\, A,\\
   u=0 &  \mbox{on}\,\, \partial A,
  \end{array}
\right.
$$
where $\varepsilon>0$, $\eta\in\mathbb{R}$, $A=\{x\in\R^{2N}:\,\,0<a<|x|<b\}$, $N\ge2$ and $f$ is a nonlinearity of $C^1$ class with supercritical growth. By a reduction argument, we show that there exists a nodal solution $u_\e$ with exactly two positive and two negative peaks, which concentrate on two different orthogonal spheres of dimension $N-1$ as $\e\rightarrow0$. In particular, we establish different concentration phenomena of four peaks when the parameter $\eta>2$, $\eta=2$ and $\eta<2$.
\end{abstract}
\maketitle
\begin{center}
\begin{minipage}{12cm}
\tableofcontents
\end{minipage}
\end{center}

\vskip0.3in
\section{Introduction}
\label{sec1}
\vskip0.1in
We study the following singularly perturbed elliptic equation with superlinear nonlinearity in an annulus in $\R^{2N}(N\ge2)$
\begin{equation}\label{eqn:qs}
\left\{
  \begin{array}{ll}
  \displaystyle  -\varepsilon^2\triangle {u}+|x|^\eta u =|x|^\eta f(u)&  \mbox{in}\,\, A,\\
   u=0 &  \mbox{on}\,\, \partial A,
  \end{array}
\right.
\end{equation}
where $\varepsilon>0, \eta\in\mathbb{R}$, $f$ is of $C^1$-class and supercritical at infinity, and $A=\{x\in\R^{2N}:\,\,0<a<|x|<b\}$. There has been plenty of results with respect to solutions with point concentration in bounded domains. Based on an energy expansion, Ni and Takagi \cite{Ni91} investigated the existence of positive solutions to the following problem with homogeneous Neumann boundary conditions and subcritical nonlinearity
\begin{equation}\label{ellip}
-\varepsilon^2\triangle {u}+u =f(u)\,\,  \mbox{in}\,\, \Omega\subset\mathbb{R}^N,\\
\end{equation}
and showed that any least energy solution has at most one local maximum, which lies on the boundary for sufficiently small $\e$. Ni and Wei \cite{Ni95} considered problem \eqref{ellip} with homogeneous Dirichlet boundary conditions and demonstrated that any least energy solution has at most one local maximal point, which concentrates around the point which stays with the maximal distance from the boundary. Dancer and Yan \cite{Dancer99} studied \eqref{ellip} for $f(s)=|s|^{p-1}s$,
$p\in\left(1, \frac{N+2}{N-2}\right)$ if $N\ge3$, $p>1$ if $N=2$. By using the Lyapunov-Schmidt reduction, the authors proved the existence of positive multi-peak solutions under the homogeneous Dirichlet boundary condition in a general domain $\Omega$ with nontrivial topology. For the further related results, we refer to \cite{Pino99,Byeon05,Aprile12,Pino20,Dancer12} and the reference therein.

In \cite{Noussair97}, Noussair and Wei studied problem \eqref{ellip} with the homogeneous Dirichlet boundary condition and $f(s)=|s|^{p-1}s$, $p\in\left(1, \frac{N+2}{N-2}\right)$ if $N\ge3$, $p>1$ if $N=2$. They obtained the existence of a least energy nodal solution and showed that the nodal solution has exactly one positive and one negative
peaks converging to two distinct points $P^1,P^2$ of $\Omega$ as $\varepsilon\rightarrow0$, respectively. D'Aprile and Pistoia \cite{Pistoia1} considered problem \eqref{ellip} with  the homogeneous Dirichlet boundary condition and established the existence of nodal solutions with multiple peaks concentrating at different points $\Omega$. We also see \cite{Alves07,Aprile11,Bartsch05,WT,Aprile13,Wu2018,Dancer11,Wei05} and the reference therein.

Ambrosetti, Malchiodi and Ni \cite{AMN1,AMN2} considered another type of concentrating solutions which concentrate on lower dimensional manifolds. They were concerned with the problem
\begin{equation}\label{eqn:eta=0}
\left\{
  \begin{array}{ll}
  \displaystyle  -\varepsilon^2\triangle {u}+V(r) u =f(u)&  \mbox{in}\,\, B,\\
   u=0 &  \mbox{on}\,\, \partial B,
  \end{array}
\right.
\end{equation}
in an annulus $B=\{x\in\R^{N}:\,\,0<a<|x|<b\}$, where $V$ is a smooth radial potential and bounded below by a positive constant.
 By introducing a modified potential $M(r)=r^{N-1}V^{\varrho}$, with $\varrho=\frac{p+1}{p-1}-\frac{1}{2}$ and $M'(b)<0$
 (respectively $M'(a)>0$), they proved there exists a family of radial positive solutions which concentrate
 on the sphere $|x|=r_\e$ with $r_\e\rightarrow b$ (respectively $r_\e\rightarrow a$) as $\e\rightarrow0$.
 They also conjectured in \cite{AMN1} that for $N\geq3$, there also exist solutions concentrating to some manifolds
 of dimension $1\leq k\leq N-2$. Concentration of positive solutions on curves was considered
 by del Pino, Kowalczyk and Wei \cite{Pino07}. It was mentioned in \cite{Ni98} that such solutions are
 of particular interest for applications to models of activator-inhibitor systems in biology.
 For related results about concentrating on higher dimensional manifolds, we can also refer to
 \cite{Malchiodi02,Malchiodi04,Bartsch08,Bartsch10} and the references therein.

By virtue of a Hopf-fibration approach, Ruf and Srikanth \cite{Ruf10} considered problem \eqref{eqn:qs}
with $\eta=0, f(s)=|s|^{p-1}s$ ($p$ maybe supercritical) in $\R^4$. They transformed the problem in an
annulus in $\R^4$ to a three-dimensional one, which can be allowed to get solutions with single point
concentration similarly to the well-known results by Ni-Wei \cite{Ni95} and del Pino-Felmer \cite{Pino99}.
Inverting the Hopf reduction, they obtained solutions concentrating on $S^1$-orbits in $\R^4$ which tend to
the inner boundary of $A$ as $\e\rightarrow0^+$. Later, Pacella and Srikanth in \cite{Pacella12} used a reduction
approach to consider \eqref{eqn:qs} with $\eta=0, f(s)=|s|^{p-1}s$, $1<p<\frac{N+3}{N-1}$ and proved the existence
of positive and sign-changing solutions concentrating on one or two $\textcolor[rgb]{1.00,0.00,0.00}{(N-1)}$ dimensional spheres.
As pointed out in \cite{Ruf10}, it seems impossible to extend the results in \cite{Ruf10} to odd dimensional cases.
Actually, only even dimensional cases are considered in \cite{Pacella12,Ruf14,Ruf14+}. In \cite{SW1},
Santra and Wei first used the Hopf fibration to study problem \eqref{ellip} with the homogeneous Dirichlet or
Neumann boundary condition in an annulus $A$ for any $N\ge2$. Actually, they reduced the problem to one in $\mathbb{R}^2$,
which yields that the nonlinearity is allowed to be polynomial growth of arbitrary order. Moreover, the solutions concentrate
on the inner boundary and outer boundary of $A$ under the Dirichlet and Neumann boundary conditions, respectively.
Via a Hopf reduction, Manna and Srikanth \cite{Manna14} considered the least energy solution of \eqref{eqn:qs}
with $f(s)=|s|^{p-1}s$ , $1<p<\frac{N+3}{N-1}$ and showed the existence of positive solutions with two peaks concentrating
along two spheres which are orthogonal to each other. For related concentrating works, we can see \cite{Manna14+,Clapp14,Ruf14+,SW2,Zhang2019}
and the references therein.

Recently, Clapp and Manna \cite{Clapp17} studied the problem
$$
-\varepsilon^{2} \Delta v+v=|v|^{p-2} v \text { in } \Omega, \quad v=0 \text { on } \partial \Omega,
$$
in domains of the form
$$
\Omega:=\left\{\left(y_{1}, y_{2}\right) \in \mathbb{R}^{n-1} \times \mathbb{R}^{m+1}:\left(y_{1},\left|y_{2}\right|\right) \in \Theta\right\},
$$
where $m \geq 1, n \geq 3$ and $\Theta$ is a bounded smooth domain in $\mathbb{R}^{n}$ with $\bar{\Theta} \subset$ $\mathbb{R}^{n-1} \times(0, \infty)$ and $p \in\left(2, \frac{2 n}{n-2}\right) .$ For particular choices of $\Theta$, the authors established the existence of single-layered and double-layered sign-changing solutions, which concentrate on spheres that converge to a single $m$-dimensional sphere contained in $\partial\Omega$. As mentioned in \cite{Clapp17}, the existence of sign-changing multi-peak solutions with more than two peaks is largely open. The aim of this paper is to exhibit some new concentration phenomena for sign-changing solutions of problem (\ref{eqn:qs}).

We state the following conditions for $f\in C^1(\R,\R)$.
\begin{itemize}
\item[ ($f_1$)]$f(t)=o(t)$ as $t\rightarrow0$, $f(t)=o(t^{p-1})$ for $p\in\left(2,\frac{2N+2}{N-1}\right)$ as $t\rightarrow\infty$, and $f$ is odd.
\item[ ($f_2$)]There exists $\mu>1$ such that $t^2f'(t)\geq\mu f(t)t>0$ for $t\not=0$.
\end{itemize}

\begin{remark}
By virtue of ($f_1$) and ($f_2$), we can conclude that
$$
tf(t)\geq(1+\mu) F(t)>0,\,\,t\not=0,\quad F(t)=\int_{0}^tf(s)ds,
$$
which is called as the well-known Ambresstti-Rabnowitz condition. Moreover, by ($f_2$) we also conclude that
$$
t\mapsto\frac{f(t)}{t} \text{\,is \,strictly\, increasing \,for\,\,} t>0,
$$
which plays an essential role in proving the existence of ground state solutions by using Nehari manifold arguments.
\end{remark}
To state the next condition, we consider the following problem in the whole space
\begin{equation}\label{eqn:spacefq}
\left\{
  \begin{array}{ll}
  \displaystyle  -\triangle {u}+u =f(u)&  \mbox{in}\,\, \R^{N+1},\\
   u>0, &  \mbox{in}\,\, \R^{N+1},
  \end{array}
\right.
\end{equation}
whose energy functional $J(u): H^{1}(\R^{N+1})\mapsto\R$ is given by
$$
J(u)=\frac{1}{2}\int_{\R^{N+1}}|\nabla u|^2+\frac{1}{2}\int_{\R^{N+1}}u^2dx
-\int_{\R^{N+1}}F(u)dx.
$$
Now we state condition ($f_3$) as follows.
\begin{itemize}
\item[ ($f_3$)] Problem (\ref{eqn:spacefq}) admits a unique positive solution $U(x)=U(|x|)$ (see \cite{Kwong89}) such that
\begin{equation}\label{eqn:Uex}
|D^{\alpha} U(x)|\leq C_1\exp(-\sigma|x|),\,x\in\R^{N+1},\text{for\,some}\,C_1,\sigma>0\,\text{and\,all}\,|\alpha|\leq 2.
\end{equation}
\end{itemize}

Let us consider $A$ under the coordinate system
$$
A=I_1\times(I_2\times S^{N-1}\times S^{N-1}),
$$
where $I_1=(a,b)$, $I_2=[0,\pi/2)$ and $S^{N-1}$ has the standard polar coordinate expression. For any $x\in A$, we write
$$
x=x(r,\theta,\theta_1^1,\theta_2^1,\cdot\cdot\cdot,\theta_{N-1}^1,\theta_1^2,\theta_2^2,\cdot\cdot\cdot,\theta_{N-1}^2),
$$
where $r\in I_1$ and $\theta\in I_2$, $\theta_1^i\in[0,2\pi)$ for $i=1,2$, $\theta_j^i\in[0,\pi)$ for $i=1,2$ and $j=2,...,N-1$.
We now state our main results.
\begin{theorem}\label{Thm:inner-out}
Assume ($f_1$)-($f_3$) hold, then problem (\ref{eqn:qs}) has a nonradial nodal solution $u_\e$ with exactly two positive and two negative peaks, which concentrate on two orthogonal spheres with dimension of $N-1$, placed by the angle $\theta=0,\theta=\frac{\pi}{2}$, belonging to the inner boundary $|x|=a$ if $\eta<2$, and belonging to the outer boundary $|x|=b$ if $\eta>2$.
\end{theorem}

\begin{theorem}\label{Thm:middle}
Assume ($f_1$)-($f_3$) hold, then problem (\ref{eqn:qs}) with $\eta=2$ has a nonradial nodal solution $u_\e$ with exactly two positive and two negative peaks. More precisely, two positive peaks concentrate on two orthogonal spheres with dimension $N-1$, placed by the angle $\theta=0,\theta=\frac{\pi}{2}$, which are contained in the surface  $|x|=\frac{1}{2}\sqrt{3a^2+b^2}$ ( or $|x|=\frac{1}{2}\sqrt{a^2+3b^2}$), and two negative peaks also concentrate on two orthogonal spheres with dimension $N-1$, placed by the angle $\theta=0,\theta=\frac{\pi}{2}$, which are contained in the surface  $|x|=\frac{1}{2}\sqrt{a^2+3b^2}$ ( or $|x|=\frac{1}{2}\sqrt{3a^2+b^2}$).
\end{theorem}

We emphasize that in \cite{Bartsch13}, Bartsch, D'Aprile and Pistoia investigated an almost critical problem with domain satisfying some certain symmetry and used a Lyapunov-Schmidt reduction scheme to construct a four-bubble nodal solution with two positive and two negative bubbles. See also \cite{Bartsch13+}. The authors in \cite{Grossi20} also studied an almost critical problem in a ball and obtained the radial nodal solutions with many bubbles concentrating at the center of the ball as $\e\rightarrow0^+$. Different from the above works, in this paper we are concerned with the supercritical case. The proof is mainly based on a Hopf-fibration reduction, energy estimates, blow-up argument, Morse index techniques and the Nehari manifold approach. Moreover, the precise locations of concentration are also considered.
\vskip0.1in
Hereafter, the letter $C$ will be repeatedly used to denote various positive constants whose exact values are irrelevant. This paper is organized as follows. Some notations and preliminary results are given in Section 2. Section 3 is devoted to the energy estimates of nodal solutions. In Section 4, we prove Theorem \ref{Thm:inner-out} and investigate the inner and outer boundary concentration. Section 5 is devoted to the proof of Theorem \ref{Thm:middle} and the concentration on the interior of the annulus $A$ is studied.

\s{Preliminary results}\label{sec:prel}
\vskip0.1in
Consider $\R^{2N}$ as the product of two copies of $\R^N$, that is, $\R^{2N}=\R^N\times\R^N$, and so we can
denote a point $x\in \R^{2N}$ by $x=(y_1,y_2)$, $y_i\in\R^{N},i=1,2$. Let us take in $\R^{N}$ spherical coordinates
$$
(\rho_1,\theta_1^1,...,\theta_{N-1}^1),\quad\quad(\rho_2,\theta_1^2,...,\theta_{N-1}^2),
$$
where $\rho_1=|y_1|$ and $\rho_2=|y_2|$ and $\theta_1^i\in[0,2\pi]$, $\theta_j^i\in[0,\pi]$ for $j=2,...,N-1$; $i=1,2$.
\begin{remark}
For $z=(z_1,...,z_N)\in\R^{N}$, we consider the spherical coordinate $(\rho,\theta_1,...,\theta_{N-1})$,
$\theta_i\in[0,\pi]$, $i=1,...,N-2$, $\theta_{N-1}\in[0,2\pi]$, $\rho=|z|$, then
\begin{equation}\label{eqn:change}
\left\{
  \begin{array}{ll}
  \displaystyle  z_1=\rho sin \theta_1...sin\theta_{N-1},\\
   z_2=\rho sin \theta_1...sin\theta_{N-2}cos\theta_{N-1},\\
   \cdot\cdot\cdot\cdot\cdot\cdot\\
   z_{N-1}=\rho sin \theta_1cos\theta_{2},\\
   z_{N}=\rho cos\theta_{1}.
  \end{array}
\right.
\end{equation}
\end{remark}
Furthermore, let us define
$$
\rho_1=r cos\theta,\quad \rho_2=r sin\theta,\quad r=|x|,\quad \theta\in[0,\frac{\pi}{2}],
$$
then for each point $x\in A$, we  have the following expression of coordinate
$$
x=(r,\theta_1^1,...,\theta_{N-1}^1,\theta_1^2,...,\theta_{N-1}^2,\theta).
$$
Note that function $u\in H_0^1(A)$ depends only on $r$ and $\theta$, that is, $u=u(r,\theta)$, and $u$ is invariant under rotations in
$y_1$, $y_2$. Thus,
$$
u\in X:=\{u\in H_0^1(A):\,\,u(x)=u(|y_1|,|y_2|)\}.
$$
Actually, the expression of the  Laplacian in these coordinates is given as
$$
\triangle_{\R^{2N}}u=u_{rr}+\frac{2N-1}{r}u_r+\frac{N-1}{r^2}u_\theta\bigg[\frac{cos\theta}{sin\theta}-\frac{sin\theta}{cos\theta}\bigg]
+\frac{u_{\theta\theta}}{r^2}.
$$

Let us define the new variables
$$
\rho=\frac{1}{2}r^2\quad\quad \varphi=2\theta
$$
and the function $v(\rho,\varphi)=u(r(\rho),\theta(\varphi))=u(\sqrt{2\rho},\frac{\varphi}{2})$. Then we have
$$
u_r=\sqrt{2\rho}v_\rho,\quad u_{rr}=2\rho v_{\rho\rho}+v_{\rho},\quad u_\theta=2v_{\varphi},\quad
u_{\theta\theta}=4v_{\varphi\varphi}.
$$
By defining $\varepsilon^2={2^{1-\frac{\eta}{2}}\varepsilon^2}$, we have $v$ satisfies
$$
-\varepsilon^2\bigg[v_{ss}+\frac{N}{s}v_s+\frac{N-1}{s^2}v_{\varphi}\frac{cos\varphi}{sin\varphi}+\frac{v_{\varphi\varphi}}{s^2}\bigg]
+\frac{v}{s^{1-\frac{\eta}{2}}}=\frac{f(v)}{s^{1-\frac{\eta}{2}}},\,\,s\in(\frac{a^2}{2},\frac{b^2}{2}),\,\varphi\in[0,\pi].
$$
So we have the reduced equation for $v$ as
\begin{equation}\label{eqn:redu-qs-}
\left\{
  \begin{array}{ll}
  \displaystyle  -\varepsilon^2\triangle {v}+\frac{v}{|x|^\alpha} =\frac{f(v)}{|x|^\alpha}&  \mbox{in}\,\, \Omega,\\
   v=0 &  \mbox{on}\,\, \partial \Omega,
  \end{array}
\right.
\end{equation}
where $\alpha=1-\frac{\eta}{2}$, $\Omega=\{x\in\R^{N+1}\big|\,\frac{a^2}{2}<|x|<\frac{b^2}{2}\}$. Obviously,  $v$ is axially symmetric.
Define $H_{\sharp}(\Omega)\subset H_0^1(\Omega)$ by
$$
H_{\sharp}(\Omega)=\{u\in H_0^1(\Omega):\,\,u(x',x_{N+1})=u(|x'|,|x_{N+1}|)\}.
$$
Observe that any nodal solution in $H_{\sharp}(\Omega)$ is axially symmetric and shall have
at leat two local maximums or minimums.
It is easy to see that $H_{\sharp}(\Omega)$ is a closed subspace of $H_0^1(\Omega)$ and
$H^1_{0,rad}(\Omega)\subset H_{\sharp}(\Omega)$. The energy functional $J_\varepsilon$ associated with (\ref{eqn:redu-qs-})
is defined as
$$
J_{\varepsilon}(u):=\frac{\varepsilon^2}{2}\int_{\Omega}|\nabla u|^2+\frac{1}{2}\int_{\Omega}\frac{u^2}{|x|^\alpha}dx
-\int_{\Omega}\frac{F(u)}{|x|^{\alpha}}dx, \quad\forall u\in H_{\sharp}(\Omega),
$$
and is of class $C^1$. The following remark implies that we can directly look for sign-changing critical points of $J_{\varepsilon}$ at
$H_{\sharp}(\Omega)$.
\begin{remark}
Let us define group $G:=O(N)\times\mathbf{Z}_2$ and the action of a topological group $G$ on a normed space
$H_0^1(\Omega)$ by a continuous map
$$
G\times H_0^1(\Omega)\rightarrow H_0^1(\Omega):\,\,g\cdot u\rightarrow gu,\quad \forall g\in G
$$
such that
$$
1\cdot u=u\quad\quad (g h)u=g(h u),\quad\quad \forall g,h\in G
$$
and
$$
u\mapsto gu \quad\text{is\,linear},\quad\quad \|gu\|_{H_0^1(\Omega)}=\|u\|_{H_0^1(\Omega)},\quad \forall g\in G,\, \text{and}\,u\in H_0^1(\Omega).
$$
Then we define space of invariant points as follows
$$
Fix(G):=\{u\in H_0^1(\Omega):\,\,gu=u, \,\forall g\in G\}.
$$
It is easy to check from the definition of $H_{\sharp}(\Omega)$ that $Fix(G)=H_{\sharp}(\Omega)$.
Hence, by the well-known principle of symmetric criticality developed by Palais (see \cite{Willem96}), we have if $u$ is a critical point of $J_{\varepsilon}$
restricted to $Fix(G)$ then $u$ is a critical point of $J_{\varepsilon}$ at $H_0^1(\Omega)$.
\end{remark}
For $d\in(\frac{a^2}{2},\frac{b^2}{2})$, let $W(x)=U(\frac{x}{d^{\alpha/2}})$ with $U\in H^1(\R^{N+1})$ defined in ($f_3$), it follows from (\ref{eqn:spacefq}) that $W$ will satisfy
\begin{equation}\label{eqn:redu-qs}
\left\{
  \begin{array}{ll}
  \displaystyle  -\triangle {u}+\frac{u}{d^\alpha} =\frac{f(u)}{d^\al}&  \mbox{in}\,\, \R^{N+1},\\
   u>0, \lim\limits_{|x|\rightarrow\infty}u(x)=0&  \mbox{in}\,\, \R^{N+1}.
  \end{array}
\right.
\end{equation}
Define the Nehari set in $H_{\sharp}(\Omega)$ corresponding to $J_\varepsilon$
as follows
$$
\mathcal{N}_\e:=\{u\in H_{\sharp}(\Omega):\,\,u^{\pm}\not\equiv0,\,\, J'_\varepsilon(u^{\pm})u^{\pm}=0\},
$$
where $u^{\pm}$ denote the positive and negative part of $u$, respectively.
We will prove the existence of the least energy nodal solutions to equation (\ref{eqn:redu-qs-}) by using the deformation technique
for $J_\e$ restricted at Nehari set $\mathcal{N}_\e$.
We can see \cite{Bartsch05+,Noussair97} for the similar results.
\begin{theorem}\label{Thm:exist} The following conclusions hold.
\begin{itemize}
\item[ (i)] For any $u,v\geq0$ belonging to $H_{\sharp}(\Omega)\setminus\{0\}$ and $uv\equiv0$,
there exist exactly two constants $t>0$ and $s>0$ such that $tu-sv\in\mathcal{N}_\e$.
\item[ (ii)]For fixed $\varepsilon>0$, equation (\ref{eqn:redu-qs-}) has a sign-changing solution
$u_\varepsilon\in H_{\sharp}(\Omega)$ such that $c_\e:=J_\varepsilon(u_\varepsilon)=\inf_{u\in\mathcal{N}_\e}J_{\varepsilon}(u)$.
\end{itemize}
\end{theorem}
\Proof  For any nonnegative functions
$u,v\in H_{\sharp}(\Omega)\setminus\{0\}$, by ($f_1$) and ($f_2$), there exist exactly two constants $t,s>0$ such that $J'_\e(tu)tu=0$ and $J'_\e(sv)sv=0$,
respectively. If the supports of $u$ and $v$ are disjoint, then $tu-sv\in\mathcal{N}_\e$.\\
It remains to prove (ii). It is easy to see from ($f_2$) that $J_\e$ is coercive and bounded from below on $\mathcal{N}_\e$.
Thus, there exists a minimizing sequence $\{u_n\}\subset \mathcal{N}_\e$ such that $J_\e(u_n)\rightarrow c_\e$ as $n\rightarrow\infty$.
It is easy to obtain from ($f_2$) that $\{u_n\}$ is bounded in $H_{\sharp}(\Omega)$. Up to subsequence, we may assume that
$$
\aligned
&u_n\rightharpoonup u_\e\,\,\text{ weakly\, in\,}\, H_{\sharp}(\Omega), \\
&u_n\rightarrow u_\e \,\,\text{strongly\, in\,}\, L^p(\Omega)\,\,\text{for}\,\,p\in (2,\frac{2N+2}{N-1}).
\endaligned
$$
 From ($f_1$)-($f_3$) we deduce that
$\int_{\Omega}|u_n^\pm|^pdx\geq C$ for some $C>0$, and so $\int_{\Omega}|u_\e^\pm|^pdx\geq C$.
Since $J'_\varepsilon(u_n^+)u_n^+=0$, by the weak lower semi-continuity,
there exists a unique $t\in(0,1]$ such that $J'_\e(tu_\e^+)tu_\e^+=0$. Similarly,
there also exists a unique $s\in(0,1]$ such that $J'_\e(su_\e^-)su_\e^-=0$. Hence,
$tu_\e^++su_\e^-\in \mathcal{N}_\e$. Based on the definition of $c_\e$, by ($f_2$) we have
$$
\aligned
c_\e&\leq J_\e(tu_\e^++su_\e^-)=J_\e(tu_\e^+)+J_\e(su_\e^-)\\
&=\frac{\mu-2}{2\mu}\bigg[t^2\int_{\Omega}(|\nabla u_\e^+|^2+|u_\e^+|^2)dx+s^2\int_{\Omega}(|\nabla u_\e^-|^2+|u_\e^-|^2)dx\bigg]\\
&\leq\frac{\mu-2}{2\mu}\bigg[\int_{\Omega}(|\nabla u_\e^+|^2+|u_\e^+|^2)dx+\int_{\Omega}(|\nabla u_\e^-|^2+|u_\e^-|^2)dx\bigg]\\
&\leq\frac{\mu-2}{2\mu}\liminf\limits_{n\rightarrow\infty}\bigg[\int_{\Omega}(|\nabla u_n|^2+|u_n|^2)dx\bigg]\\
&=c_\e,
\endaligned
$$
which implies that $t=s=1$. That is to say, $u_\e\in \mathcal{N}_\e$ and $J_\e(u_\e)=c_\e$. It suffices to prove
$J'_\e(u_\e)\varphi=0$ for any $\varphi\in C_0^\infty(\Omega)$. Assume by contradiction that $J'_\e(u_\e)\not=0$, then by the
continuity of $J'_\e$, there exist
$\delta>0$ and $c>0$ such that $\|J'_\e(v)\|\geq c$
if $\|v-u_\e\|\leq 3\delta$. Define $D:=(\frac{1}{2},\frac{3}{2})\times(\frac{1}{2},\frac{3}{2})$ and
$l(s,t)=su_\e^++tu_\e^-$ for $(s,t)\in D$. Note that, by ($f_2$), we have for $(t,s)\not=(1,1)$
\begin{equation}\label{eqn:exist1}
J_\e(su_\e^++tu_\e^-)=J_\e(su_\e^+)+J_\e(tu_\e^-)<J_\e(u_\e^+)+J_\e(u_\e^-)=c_\e,
\end{equation}
which yields that
$$
\bar{c}:=\max\limits_{(t,s)\in\partial D}J_\e\circ l(s,t)<c_\e.
$$
Let us set $\nu:=\min\{\frac{c_\e-\bar{c}}{2},\frac{c\delta}{8}\}$ and denote $B(u_\e,\delta)$ by the ball in $H_{\sharp}(\Omega)$
of radius
$\delta$ centered at $u_\e$. Arguing as Lemma 2.3 in \cite{Willem96}, we obtain a deformation $\gamma$ satisfying
\begin{itemize}
\item[ (a)] $\gamma(1,u)=u$ if $u\not\in J_\e^{-1}([c_\e-2\nu,c_\e+2\nu])$,
\item[ (b)] $\gamma(1,J_\e^{c_\e+\nu}\cap B(u_\e,\delta))\subset J_\e^{c_\e-\nu}$,
\item[ (c)] $J_\e(\gamma(1,u))\leq J_\e(u)$ for all $u\in H_{\sharp}(\Omega)$.
\end{itemize}
Based on the above facts, we have immediately
\begin{equation}\label{eqn:exist2}
\max\limits_{(t,s)\in \bar{D}}J_\e(\gamma(1,l(t,s)))<c_\e.
\end{equation}
Let us define $g(s,t):=\gamma(1,l(s,t))$ and
$$
\aligned
&G_0(s,t):=\bigg(J'_\e(t u_\e^+)u_\e^+,J'_\e(s u_\e^-)u_\e^-\bigg),\\
&G_1(s,t):=\bigg(\frac{1}{t}J'_\e(g^+(s,t))g^+(s,t),\frac{1}{s}J'_\e(g^-(s,t))g^-(s,t)\bigg).
\endaligned
$$
Since $u_\e\in \mathcal{N}_\e$, by $(f_1)$ and $(f_2)$, we have $\deg(G_0, D,\mathbf{0})=1$.
By virtue of the definition of $\bar{c}$ and conclusion (a), we have $ g\equiv l$ on $\partial D$.
And so, we immediately obtain $\deg(G_1, D,\mathbf{0})=\deg(G_0, D,\mathbf{0})=1$, It yields $G_1(s,t)=0$ for some
$(s,t)\in D$. So, $g(s,t):=\gamma(1,l(s,t))\in \mathcal{N}_\e$ which yields a contradiction
by combining (\ref{eqn:exist2}) and the definition of $c_\e$. The proof is complete.
\qed

In order to prove that sign-changing solution $u_\e\in H_{\sharp}(\Omega)$ of equation (\ref{eqn:redu-qs-})
in Theorem \ref{Thm:exist} is a nonradial symmetric solution, we need to estimate the Morse index
with respect to energy functional $J_\e$ which is defined at $H_{\sharp}(\Omega)$.
Consider the Hilbert space $H:=H_{\sharp}(\Omega)\cap H^2(\Omega)$, which is endowed with the scalar product from
$H^2(\Omega)$. Denote by $\|\cdot\|_{H}$ the induced norm. Define functions
$$
\aligned
&\Upsilon_{\pm}: H\rightarrow\R,\quad \Upsilon_{\pm}(u):=\int_{\Omega}|\nabla u^{\pm}|^2dx=
\int_{\Omega}\nabla u\cdot\nabla u^{\pm}dx,\\
&\Psi_{\pm}: H\rightarrow\R,\quad \Psi_{\pm}(u):=\int_{\Omega}\frac{1}{|x|^\al}(-u+f(u))u^\pm dx.
\endaligned
$$
We recall some results in \cite{Bartsch03} where some properties of the
above functionals were obtained at space $H_0^1(\Omega)\cap H^2(\Omega)$. The same proof is valid for the following Lemma.
\begin{lemma}\label{Lem:mox} The following conclusions hold.
\begin{itemize}
\item[ (i)] $\Upsilon_{\pm}\in C^1(H)$  with derivative $\Upsilon'_{\pm}(u)\in H^{-1}_{\sharp}(\Omega)$ given by
$$
\Upsilon'_{\pm}(u)v:=
\int_{\pm u>0}((-\triangle u) v+\nabla u\nabla v)dx.
$$
\item[ (ii)]  $\Psi_{\pm}\in C^1(H)$ with derivative given by
$$
\Psi'_{\pm}(u)v:=\int_{\Omega}\frac{1}{|x|^\al}[(-1+f'(u^\pm))u^\pm v+(-u^\pm+f(u^\pm))v] dx.
$$
\item[ (iii)]  The set $\mathcal{N}_\e\cap H$ is a $C^1$-manifold of codimension $2$ in $H$.
\item[ (iv)]  $m(u_\e)$=2, where $u_\e$ is a sign changing solution obtained in Theorem \ref{Thm:exist} and
$m(u)$ denotes the Morse index of critical point $u$ of $J_\e$ at $H_{\sharp}(\Omega)$.
\end{itemize}
\end{lemma}

\begin{lemma}\label{Lem:nonradial} Nodal solution $u_\e\in H_{\sharp}(\Omega)$ of equation (\ref{eqn:redu-qs-}) is nonradial.
\end{lemma}
\Proof
Let $u\in H_{0,rad}(\Omega)$ be a radial nodal solution of (\ref{eqn:redu-qs-}). Then by the elliptic regularity estimates, we have
$u\in H$. Using conclusions (i) and (ii) of Lemma \ref{Lem:mox} and ($f_2$), a direct computation yields
\begin{equation}\label{eqn:non1}
\aligned
J_\e''(u)(u^{\pm},u^{\pm})&=\int_{\Omega}\e^2|\nabla u^{\pm}|^2dx
+\int_{\Omega}\frac{(u^{\pm})^2-f'(u^{\pm})(u^{\pm})^2}{|x|^2}dx\\
&\leq-(\mu-1)\int_{\Omega}\frac{f(u^{\pm})u^{\pm}}{|x|^2}dx<0.
\endaligned
\end{equation}
On the other hand, it follows form Lemma 2.2 in \cite{Manna14} that we can construct a $v\in H_{\sharp}(\Omega)$ by taking
$v(r,\varphi)=u(r)(c+cos(2\varphi))$ for $\varphi\in[0,\pi]$,
where the co-ordinate system of
$\Omega$ is taken as the standard polar co-ordinate system and $c$ is chosen so that $u$ and $v$ are
orthogonal in $L^2(\Omega)$. That is to say,
$$
\int_{a^2/2}^{b^2/2}\int_{0}^{\pi}u^2(r)(c+cos(2\varphi))r^N sin^{N-1}\varphi drd\varphi=0.
$$
$\bigtriangleup v$ takes in the standard polar coordinate system the form
$$
\triangle v=v_{rr}+\frac{N}{r}v_r+\frac{1}{r^2}v_{\varphi\varphi}+\frac{N-1}{r^2}\frac{cos\varphi}{sin\varphi}v_\varphi,
$$
and then
$$
\aligned
-&\e^2\triangle v+\frac{1}{|x|^\al}(1-f'(u))v\\
&=\frac{1}{|x|^\al}(f(u)-f'(u)u)(c+cos2\varphi)
+\frac{4\e^2}{r^2}u\,cos2\varphi+\frac{4(N-1)\e^2}{r^2}u\,cos^2\varphi.
\endaligned
$$
By making the change of variables and integrating by parts, we have
\begin{equation}\label{eqn:non2}
\aligned
J_\e''(u)(v,v)=&\int_{\Omega}\frac{1}{r^\al}[f(u)u-f'(u)u^2](c+cos2\varphi)^2dx
+\int_{\Omega}\frac{4\e^2}{r^2}u^2cos2\varphi(c+cos2\varphi)dx\\
&+\int_{\Omega}\frac{4(N-1)\e^2}{r^2}u^2cos^2\varphi(c+cos2\varphi)dx\\
=&C_1\int_{a^2/2}^{b^2/2}[f(u)u-f'(u)u^2]r^{N-\al}dr
+4\e^2 C_2\int_{a^2/2}^{b^2/2}r^{N-2}u^2dr\\
&+4(N-1)\e^2 C_3\int_{a^2/2}^{b^2/2}r^{N-2}u^2dr
\endaligned
\end{equation}
for some $C_1, C_2,C_3>0$ independently of $\e$. Take $C_4=(a^2/2)^{\al-2}$ if $\al<2$;
$C_4=1$ if $\al=2$ and $C_4=(b^2/2)^{\al-2}$ if $\al>2$.
Then by (\ref{eqn:non2}) and ($f_2$), we have
$$
\aligned
J_\e''(u)(v,v)&\leq C_1\int_{a^2/2}^{b^2/2}[f(u)u-f'(u)u^2]r^{N-\al}dr
+4\e^2 [C_2+(N-1)C_3]C_4\int_{a^2/2}^{b^2/2}r^{N-\al}u^2dr\\
&\leq C_1\int_{a^2/2}^{b^2/2}(1-\mu)f(u)ur^{N-\al}dr
+4\e^2 [C_2+(N-1)C_3]C_4\int_{a^2/2}^{b^2/2}r^{N-\al}u^2dr\\
&\leq C_1(\mu-1)\int_{a^2/2}^{b^2/2}[u^2-f(u)u]r^{N-\al}dr\\
&= -C(\mu-1)\e^2\int_{\Omega}|\nabla u|^2dx<0
\endaligned
$$
for $C>0$ independently of $\e$, where the third inequality
we choose $\e$ small enough so that
$$
4\e^2 [C_2+(N-1)C_3]C_4<C_1(\mu-1).
$$
We now prove that $\{u^+,u^-,v\}\subset H_{\sharp}(\Omega)$ are linearly independent vectors.
If $v$ and $u^+$ are linearly dependent, then $v=k u^+$ for $k\in\R$. Clearly, $k\not=0$.
According to $u$ and $v$ being orthogonal, we have
$$
\int_{\Omega}uv dx=\int_{\Omega}k u u^+ dx=k\int_{\Omega}|u^+|^2dx>0.
$$
It is a contradiction. Similarly, we can obtain $v, u^-$ are linearly independent. Hence,
we see that the Morse index of any radial nodal solution
$u$ is greater or equal to $3$. Recalling Lemma \ref{Lem:mox} (iv), we deduce that $u_\e\in H_{\sharp}(\Omega)$ is nonradial.
The proof is complete.
\qed

\s{The energy estimates}
\vskip0.1in
We now state an upper estimate for energy $c_\e$ defined in Theorem \ref{Thm:exist}.
\begin{lemma}\label{Lem:upper} For small $\e>0$, there exist $d_1,d_2>0$ such that
$$
c_\e \leq2\e^{N+1}[(d_1^{(N-1)\frac{\al}{2}}+d_2^{(N-1)\frac{\al}{2}})J(U)+O(\e)],
$$
where $U\in H^1(\R^{N+1})$ is the unique positive solution of equation (\ref{eqn:spacefq}).
\end{lemma}
\Proof
Define
$$
\aligned
&\mathcal{A}_{\varepsilon}(x)=\phi(\frac{x-P}{d_1^{\al/2}})W_1(\frac{x-P}{\e })+\phi(\frac{x+P}{d_1^{\al/2}})W_1(\frac{x+P}{\e }),\\
&\mathcal{B}_{\varepsilon}(x)=\phi(\frac{x-Q}{d_2^{\al/2}})W_2(\frac{x-Q}{\e })+\phi(\frac{x+Q}{d_2^{\al/2}})W_2(\frac{x+Q}{\e }),
\endaligned
$$
where $P=(\mathbf{0},d_1)$ and $Q=(\mathbf{0},d_2)$ for some $d_1,d_2\in(a^2/2,b^2/2)$ and $d_1\not=d_2$,
and
 $\phi$ is a non-negative
smooth radial function supported in $B_{2\gamma}(0)$ with $|\nabla \phi|\leq \frac{2}{\gamma}$ and
$$
\phi(r)=\left\{
  \begin{array}{ll}
  \displaystyle  1&  \mbox{for}\,\, r\in[0,\gamma],\\
   0,&  \mbox{for}\,\, r\in[2\gamma,+\infty),
  \end{array}
\right.
$$
where $\gamma$ is chosen so that
$$
\max\bigg\{2\gamma d_1^{\frac{\al}{2}},2\gamma d_2^{\frac{\al}{2}}\bigg\}
<\min\{dist(P,\partial\Omega),dist(Q,\partial\Omega)),\frac{|d_1-d_2|}{2}\}.
$$
Note that $W_1$ and $W_2$ are least energy positive solutions of equation (\ref{eqn:redu-qs}) with $d=d_1$ and $d=d_2$, respectively. Clearly, $\mathcal{A}_{\varepsilon}, \mathcal{B}_{\varepsilon}\in H_{\sharp}(\Omega)$.
Observe that the supports of $\phi(\frac{x-P}{d_1^{\al/2}})W_1(\frac{x-P}{\e })$, $\phi(\frac{x-P}{d_1^{\al/2}})W_1(\frac{x+P}{\e })$,
$\phi(\frac{x-Q}{d_2^{\al/2}})W_2(\frac{x-Q}{\e })$ and $\phi(\frac{x+Q}{d_2^{\al/2}})W_2(\frac{x+Q}{\e })$ are disjoint each other
for small $\e$. Based on the definition of $\mathcal{N}_\e$, there exist exactly two constants $t_\e,s_\e>0$ such that
\begin{equation}\label{eqn:defineC}
\mathcal{C}_{\varepsilon}:=t_\e\mathcal{A}_{\varepsilon}-s_\e\mathcal{B}_{\varepsilon}\in \mathcal{N}_\e,
\end{equation}
whose energy is given as follows
$$
J_{\varepsilon}(\mathcal{C}_{\varepsilon})=\frac{\varepsilon^2}{2}\int_{\Omega}|\nabla \mathcal{C}_{\varepsilon}|^2dx+\frac{1}{2}\int_{\Omega}\frac{|\mathcal{C}_{\varepsilon}|^2}{|x|^\alpha}dx
-\int_{\Omega}\frac{F(\mathcal{C}_{\varepsilon})}{|x|^{\alpha}}dx.
$$
According to (\ref{eqn:defineC}), one has
\begin{equation}\label{eqn:hahri2}
\aligned
&\frac{\varepsilon^{1-N}}{2}\int_{\Omega}|\nabla \mathcal{C}_{\varepsilon}|^2dx\\
=&
\frac{\varepsilon^{1-N}t_\e^2}{2}\int_{\Omega}|\nabla \mathcal{A}_{\varepsilon}|^2dx
+\frac{\varepsilon^{1-N}s_\e^2}{2}\int_{\Omega}|\nabla \mathcal{B}_{\varepsilon}|^2dx\\
=&\frac{\varepsilon^{1-N}t_\e^2}{2}\int_{\Omega}\bigg|\nabla \bigg(\phi(\frac{x-P}{d_1^{\al/2}})W_1(\frac{x-P}{\e })\bigg)\bigg|^2dx+
\frac{\varepsilon^{1-N}t_\e^2}{2}\int_{\Omega}\bigg|\nabla \bigg(\phi(\frac{x+P}{d_1^{\al/2}})W_1(\frac{x+P}{\e })\bigg)\bigg|^2dx\\
&+\frac{\varepsilon^{1-N}s_\e^2}{2}\int_{\Omega}\bigg|\nabla \bigg(\phi(\frac{x-Q}{d_2^{\al/2}})W_2(\frac{x-Q}{\e })\bigg)\bigg|^2dx+
\frac{\varepsilon^{1-N}s_\e^2}{2}\int_{\Omega}\bigg|\nabla \bigg(\phi(\frac{x+Q}{d_2^{\al/2}})W_2(\frac{x+Q}{\e })\bigg)\bigg|^2dx\\
=&:I_1+I_2+I_3+I_4.
\endaligned
\end{equation}
A direct computation yields
$$
\aligned
I_1=I_2&=\frac{1}{2}d_1^{(N-1)\frac{\al}{2}}t_\e^{2}\int_{B_{2\gamma/\e}(0)}\bigg|\nabla \bigg(\phi(\e y)U(y)\bigg)\bigg|^2dy\\
&=\frac{1}{2}d_1^{(N-1)\frac{\al}{2}}t_\e^{2}\bigg[\int_{\R^{N+1}}|\nabla U|^2dy+O(\e)\bigg],\\
\endaligned
$$
and
$$
\aligned
I_3=I_4&=\frac{1}{2}d_1^{(N-1)\frac{\al}{2}}s_\e^{2}\int_{B_{2\gamma/\e}(0)}\bigg|\nabla \bigg(\phi(\e y)U(y)\bigg)\bigg|^2dy\\
&=\frac{1}{2}d_2^{(N-1)\frac{\al}{2}}s_\e^{2}\left[\int_{\R^{N+1}}|\nabla U|^2dy+O(\e)\right],
\endaligned
$$
where we use the fact that $W_1(\cdot)=U(\frac{\cdot}{d_1^{\al/2}})$ and $W_2(\cdot)=U(\frac{\cdot}{d_2^{\al/2}})$.
Using the dominate convergence theorem, the following estimates hold true
\begin{equation}\label{eqn:hahri3}
\aligned
&\e^{-(1+N)}\int_{\Omega}\frac{\mathcal{C}_{\varepsilon}^2}{2|x|^\alpha}dx\\
=&
\e^{-(1+N)}t^2_\e\int_{\Omega}\frac{|\mathcal{A}_{\varepsilon}|^2}{2|x|^\al}dx
+\e^{-(1+N)}s^2_\e\int_{\Omega}\frac{|\mathcal{B}_{\varepsilon}|^2}{2|x|^\al}dx\\
=&\e^{-(1+N)}t^2_\e\int_{\Omega}\frac{|\phi(\frac{x-P}{d_1^{\al/2}})W_1(\frac{x-P}{\e })|^2}{2|x|^\al}dx
+\e^{-(1+N)}t^2_\e\int_{\Omega}\frac{|\phi(\frac{x+P}{d_1^{\al/2}})W_1(\frac{x+P}{\e })|^2}{2|x|^\al}dx\\
&+\e^{-(1+N)}s^2_\e\int_{\Omega}\frac{|\phi(\frac{x-Q}{d_2^{\al/2}})W_2(\frac{x-Q}{\e })|^2}{2|x|^\al}dx+
\e^{-(1+N)}s^2_\e\int_{\Omega}\frac{|\phi(\frac{x+Q}{d_2^{\al/2}})W_2(\frac{x+Q}{\e })|^2}{2|x|^\al}dx\\
=&d_1^{(N-1)\frac{\al}{2}}t^2_\e\bigg[\int_{B_{2\gamma/\e}(0)}\frac{|\phi(\e y)U(y)|^2}{2|P+\e d_1^{\al/2}y|^\al}dy
+\int_{B_{2\gamma/\e}(0)}\frac{|\phi(\e y)W_1(\frac{x+P}{\e })|^2}{2|\e d_1^{\al/2}y-P|^\al}dy\bigg]\\
&+d_2^{(N-1)\frac{\al}{2}}s^2_\e\bigg[\int_{B_{2\gamma/\e}(0)}\frac{|\phi(\e y)U(y)|^2}{2|Q+\e d_2^{\al/2}y|^\al}dy+
\int_{B_{2\gamma/\e}(0)}\frac{|\phi(\e y)U(y)|^2}{2|\e d_2^{\al/2}y-Q|^\al}dy\bigg]\\
=&d_1^{(N-1)\frac{\al}{2}}\frac{t_\e^{2}}{2}\bigg[\int_{\R^{N+1}}U^2dx+O(\e)\bigg]+
d_2^{(N-1)\frac{\al}{2}}\frac{s_\e^2}{2}\bigg[\int_{\R^{N+1}}U^2dx+O(\e)\bigg].
\endaligned
\end{equation}
Moreover, we have
\begin{equation}\label{eqn:hahri4}
\aligned
&\e^{-(1+N)}\int_{\Omega}\frac{F(\mathcal{C}_{\varepsilon})}{|x|^{\alpha}}dx\\
=&
\e^{-(1+N)}\int_{\Omega}\frac{F(\mathcal{A}_{\varepsilon})}{|x|^\al}dx
+\e^{-(1+N)}\int_{\Omega}\frac{F(\mathcal{B}_{\varepsilon})}{|x|^\al}dx\\
=&\e^{-(1+N)}\int_{\Omega}\frac{F(t_\e\phi(\frac{x-P}{d_1^{\al/2}})W_1(\frac{x-P}{\e }))}{|x|^\al}dx
+\e^{-(1+N)}\int_{\Omega}\frac{F(t_\e\phi(\frac{x+P}{d_1^{\al/2}})W_1(\frac{x+P}{\e }))}{|x|^\al}dx\\
&+\e^{-(1+N)}\int_{\Omega}\frac{F(s_\e\phi(\frac{x-Q}{d_2^{\al/2}})W_2(\frac{x-Q}{\e }))}{|x|^\al}dx+
\e^{-(1+N)}\int_{\Omega}\frac{F(s_\e\phi(\frac{x+Q}{d_2^{\al/2}})W_2(\frac{x+Q}{\e }))}{|x|^\al}dx.
\endaligned
\end{equation}
By Fatou's lemma and ($f_2$), if $t_\e\rightarrow\infty$ as $\e\rightarrow0$, one has
$$
\aligned
\frac{\e^{-(1+N)}}{t_\e^2}\int_{\Omega}\frac{F(t_\e\phi(\frac{x-P}{d_1^{\al/2}})W_1(\frac{x-P}{\e }))}{|x|^\al}dx
=d_1^{(N-1)\frac{\al}{2}}\int_{B_{2\gamma/\e}(0)}
\frac{F(t_\e\phi(\e y)U(y))}{t_\e^2|P+\e d_1^{\frac{\al}{2}}y|^\al}dy\rightarrow\infty,
\endaligned
$$
and if $s_\e\rightarrow\infty$ as $\e\rightarrow0$,
$$
\aligned
\frac{\e^{-(1+N)}}{s_\e^2}\int_{\Omega}\frac{F(s_\e\phi(\frac{x-Q}{d_2^{\al/2}})W_2(\frac{x-Q}{\e }))}{|x|^\al}dx
=d_2^{(N-1)\frac{\al}{2}}\int_{B_{2\gamma/\e}(0)}
\frac{F(s_\e\phi(\e y)U(y))}{s_\e^2|Q+\e d_2^{\frac{\al}{2}}y|^\al}dy\rightarrow\infty.
\endaligned
$$
Then it follows from \eqref{eqn:hahri2} that $t_\e,s_\e$ are bounded uniformly for $\e$. Assume that
\begin{equation}\label{eqn:epsilon}
t_\e\rightarrow t_0\geq0,\quad\quad s_\e\rightarrow s_0\geq0,
\end{equation}
as $\e\rightarrow0^+$. By the dominate convergence theorem,
$$
\aligned
\e^{-(1+N)}\int_{\Omega}\frac{F(t_\e\phi(\frac{x-P}{d_1^{\al/2}})W_1(\frac{x-P}{\e }))}{|x|^\al}dx
&=d_1^{(N-1)\frac{\al}{2}}\int_{B_{2\gamma/\e}(0)}
\frac{F(t_\e\phi(\e y)U(y))}{|P+\e d_1^{\frac{\al}{2}}y|^\al}dy\\
&=d_1^{(N-1)\frac{\al}{2}}\bigg[\int_{\R^{N+1}}
F(t_0U(y))dy+O(\e)\bigg].
\endaligned
$$
Similarly, we also have
\begin{equation}\label{eqn:nahar5}
\aligned
&\e^{-(1+N)}\int_{\Omega}\frac{F(t_\e\phi(\frac{x+P}{d_1^{\al/2}})W_1(\frac{x+P}{\e }))}{|x|^\al}dx
=d_1^{(N-1)\frac{\al}{2}}\bigg[\int_{\R^{N+1}}
F(t_0 U(y))dy+O(\e)\bigg],\\
&\e^{-(1+N)}\int_{\Omega}\frac{F(s_\e\phi(\frac{x-Q}{d_2^{\al/2}})W_2(\frac{x-Q}{\e }))}{|x|^\al}dx
=d_2^{(N-1)\frac{\al}{2}}\bigg[\int_{\R^{N+1}}
F(s_0 U(y))dy+O(\e)\bigg],\\
&\e^{-(1+N)}\int_{\Omega}\frac{F(s_\e\phi(\frac{x+Q}{d_2^{\al/2}})W_2(\frac{x+Q}{\e }))}{|x|^\al}dx
=d_2^{(N-1)\frac{\al}{2}}\bigg[\int_{\R^{N+1}}
F(s_0 U(y))dy+O(\e)\bigg].
\endaligned
\end{equation}
In view of the above facts, we immediately have
\begin{equation}\label{eqn:nahar6}
\aligned
&\e^{-(N+1)}J_{\e}(\mathcal{C}_{\varepsilon})=2d_1^{(N-1)\frac{\al}{2}}\bigg[\frac{t_0^{2}}{2}\int_{\R^{N+1}}
|\nabla U|^2dx+\frac{t_0^{2}}{2}\int_{\R^{N+1}}
U^2dx-\int_{\R^{N+1}}
F(t_0 U)dx+O(\e)\bigg]\\
&+2d_2^{(N-1)\frac{\al}{2}}\bigg[\frac{s_0^{2}}{2}\int_{\R^{N+1}}
|\nabla U|^2dx+\frac{s_0^{2}}{2}\int_{\R^{N+1}}
U^2dx-\int_{\R^{N+1}}
F(s_0 U)dx+O(\e)\bigg].
\endaligned
\end{equation}
Based on the definition of $c_\e$ and recalling Theorem \ref{Thm:exist}, we have
\begin{equation}\label{eqn:nahar7}
\aligned
&\e^{-(N+1)}c_\e\\
=&2d_1^{(N-1)\frac{\al}{2}}\max\limits_{t\in(0,+\infty)}\bigg[\frac{t^{2}}{2}\int_{\R^{N+1}}
|\nabla U|^2dx+\frac{t^{2}}{2}\int_{\R^{N+1}}
U^2dx-\int_{\R^{N+1}}
F(t U)dx+O(\e)\bigg]\\
&+2d_2^{(N-1)\frac{\al}{2}}\max\limits_{s\in(0,+\infty)}\bigg[\frac{s^{2}}{2}\int_{\R^{N+1}}
|\nabla U|^2dx+\frac{s^{2}}{2}\int_{\R^{N+1}}
U^2dx-\int_{\R^{N+1}}
F(s U)dx+O(\e)\bigg]\\
\leq&2[d_1^{(N-1)\frac{\al}{2}}+d_2^{(N-1)\frac{\al}{2}}]J(U)+O(\e).
\endaligned
\end{equation}
The proof is complete.
\qed

Since $u_\e\in\mathcal{N}_\e$ and $u_\e^{\pm}\not\equiv0$, by Sobolev's inequality and ($f_1$),
we have for any $\epsilon>0$, there exists $C_\epsilon>0$ such that
$$
\aligned
C\|u_\e^{\pm}\|_p^{\frac{2}{p}}&\leq \int_{\Omega}|\nabla u_\e^{\pm}|^2+\frac{|u_\e^{\pm}|^2}{|x|^\al}dx
\leq \int_{\Omega}\frac{F(u_\e^{\pm})}{|x|^\al}dx\\
&\leq\frac{2}{a^\al}\int_{\Omega}(\epsilon|u_\e^{\pm}|^2+C_\epsilon|u_\e^{\pm}|^p)dx,
\endaligned
$$
from which we deduce that there exists $C>0$ independent of $\e$ such that
$$
\int_{\Omega}|u_\e^{\pm}|^pdx\geq C>0.
$$
By Sobolev's embedding, Moser's iteration,
 we can show $u_\e\in L^\infty(\Omega)$, and then by the elliptic estimates, $u_\e\in C^1(\bar{\Omega})$. Indeed, it follows from (\ref{eqn:nahar7})
 that there exists $C>0$ independent of $\e$(small enough) such that $\max_{x\in\Omega}|u_\e(x)|\leq C$.
Let $u_\e(P_\e)=\max_{x\in\Omega}u_\e(x)$ and $u_\e(Q_\e)=\min_{x\in\Omega}u_\e(x)$ for
some $P_\e,Q_\e\in\Omega$.

Define
\begin{equation}\label{eqn:omega+}
\aligned
&\Omega^1:=\{x=(\mathbf{x},x_{N+1})\in\Omega:\,\,x_{N+1}>0\}\textcolor[rgb]{1.00,0.00,0.00}{,}\\
&\Omega^2:=\{x=(\mathbf{x},x_{N+1})\in\Omega:\,\,x_{N+1}<0\}\textcolor[rgb]{1.00,0.00,0.00}{,}\\
&\Omega^o:=\{x=(\mathbf{x},x_{N+1})\in\Omega:\,\, x_{N+1}=0\},
\endaligned
\end{equation}
then $\Omega=\Omega^1\cup\Omega^2\cup\Omega^o$. Assume that
\begin{equation}\label{eqn:limit}
\lim\limits_{\e\rightarrow 0^+}P_\e= \bar{P},\quad\quad \lim\limits_{\e\rightarrow 0^+}Q_\e= \bar{Q}
\end{equation}
for $\bar{P},\bar{Q}\subset\bar{\Omega}$.
\begin{lemma}\label{Lem:axis}
As $\e\rightarrow0^+$, we have
\begin{itemize}
\item[ (i)]$
\frac{\min\{dist(\pm P_\e,\partial\Omega), dist(\pm Q_\e,\partial\Omega)\}}{\e}\rightarrow\infty,
$ and
\item[ (ii)]$
\frac{\min\{dist(\pm P_\e,\Omega^o), dist(\pm Q_\e,\Omega^o)\}}{\e}\rightarrow\infty.
$
\end{itemize}
\end{lemma}
\Proof
We firstly prove the first conclusion of Lemma \ref{Lem:axis}. We use similar arguments as Theorem 2.2 in \cite{Ni95}
to show conclusion (i), see also \cite{Ruf10}.
We firstly prove $dist\{P_\e,\partial\Omega\}/\e\rightarrow\infty$ as $\e\rightarrow0$.
Since ${\e}^2\triangle u_{\e}=\frac{1}{|x|^\al}(u_{\e}-f(u_{\e}))$ and ${\e}^2\triangle u_{\e}(P_{\e})\leq0$, we have
$u_{\e}(P_{\e})\leq f(u_{\e}(P_{\e}))$ which implies by conditions ($f_1$) that there exists $c_1>0$ such that
$u_{\e}(P_{\e})\geq c_1$.
Assume on the contrary that there exists $c>0$ such that $\lim_{\e\rightarrow0}dist\{P_\e,\partial\Omega\}/\e<c$.
Thus, we may assume that $P_\e\rightarrow \bar{P}\in \partial\Omega$,
that is, $|\bar{P}|=a^2/2$ or $|\bar{P}|=b^2/2$.
Using a ``boundary straightening" around the point $\bar{P}$, we may assume that $\bar{P}$ is
the origin and the inner normal to $\partial\Omega$ at $\bar{P}$
is the direction of the positive $x_{N+1}$-axis.
Define $w_\e(x)=u_\e(\mathcal{G}(q_\e+\e x))$ where $\mathcal{G}$ is a ``straightening map''
defined by
$$
\mathcal{G}: \bar{B}_{\kappa/\e}(0)\cap \{x_{N+1}\geq-\alpha_\e\}\subset\R^{N+1}\rightarrow \Omega
$$
and
$\mathcal{G}(q_\e)=P_\e$. Here, $\kappa>0$ and $\alpha_\e>0$ is bounded and $\al_\e\rightarrow\al\geq0$. By elliptic regularity theory,
$w_\e\rightarrow w_0$ in $C_{loc}^2(\R^{N+1}_{\al,+})$ (see \cite{Ni95}), where
$\R^{N+1}_{\al,+}=\{x(\mathbf{x},x_{N+1})\in\R^{N+1}|\,x_{N+1}>-\al\}$, and
$w_0$ satisfies
$$
\left\{
  \begin{array}{ll}
  \displaystyle  -\triangle {w}+\frac{w}{|\bar{P}|^\alpha} =\frac{f(w)}{|\bar{P}|^\alpha}, \,\, w>0, &\mbox{in}\,\,\R^{N+1}_{\al,+},\\
   w(x)=0,\,\,& \mbox{on}\,\,\partial\R^{N+1}_{\al,+}.
  \end{array}
\right.
$$
Form Theorem 1.1 in \cite{Lions82}, we deduce that $w_0\equiv0$. It contradicts with
$w_0(0)=\lim_{\e\rightarrow 0}w_\e(0)=\lim_{\e\rightarrow 0}u_\e(P_\e)\geq c_1$.
The other cases can be proved similarly.

 It remains to prove the second conclusion of Lemma \ref{Lem:axis}. We first claim $P_\e,Q_\e\not\in \Omega^o$ for $\e$ small enough.
 Here we only prove that $P_\e\not\in \Omega^o$ for $\e$ small enough, and similar arguments to $Q_\e$ hold true.
Assume by contradiction that there exists $\{\e_n\}$ satisfying $\e_n\rightarrow 0^+$ such that
\begin{equation}\label{eqn:limit0}
P_{\e_n}=(\mathbf{x}_{\e_n},0)\in\Omega^o
\end{equation}
with $|\mathbf{x}_{\e_n}|=d_{\e_n}$, and such that $d_{\e_n}\rightarrow d$ for some $d\geq\frac{a^2}{2}$, as $n\rightarrow +\infty$.
Clearly, $u_{\e_n}(P_{\e_n})\geq c_1$ and $|\mathbf{x}_{\e_n}|>\frac{d}{2}$ for large $n$.
Thus, for fixed $n$, we can always find $k$ points
$$
P^i_{\e_n}=(\mathbf{x}_{\e_n}^i,0)\in \Omega\quad\text{with}\,\, |\mathbf{x}_{\e_n}^i|=|\mathbf{x}_{\e_n}|,\,\,i=1,...,k
$$
such that
$$
|P^i_{\e_n}-P^j_{\e_n}|\geq c_0,\,\,\, i\not=j,\,\,i,j=1,...,k
$$
for some $c_0>0$. Clearly, $|P^i_{\e_n}|=|P_{\e_n}|$. The definition of $H_{\sharp}(\Omega)$ implies that $u_{\e_n}(P_{\e_n})=u_{\e_n}(P^i_{\e_n})$, and so
$u_{\e_n}(P^i_{\e_n})\geq c_1$ for $i=1,...,k$. Recalling Theorem \ref{Thm:exist}, we have for any given $R>0$
\begin{equation}\label{eqn:zhu1}
\aligned
c_{\e_n}&=\int_{\Omega}\frac{1}{|x|^\al}(\frac{1}{2}f(u_{\e_n})u_{\e_n}-F(u_{\e_n}))dx\\
&\geq\sum_{i=1}^{k}\int_{B_{{\e_n} R}(P^i_{\e_n})}\frac{1}{|x|^\al}(\frac{1}{2}f(u_{\e_n})u_{\e_n}-F(u_{\e_n}))dx.
\endaligned
\end{equation}
Here, we have used the fact that the supports of $B_{{\e_n} R}(P^i_{\e_n}), i=1,...,k$ are disjoint
for any given $R>0$.
By making the changes of variables $v^i_{\e_n}(y)=u_{\e_n}({\e_n} y+P^i_{\e_n})$,
using elliptic regularity theory and the boundedness of $\{u_{\e_n}\}$ in $H_{\sharp}(\Omega)$,
we obtain
$v^i_{\e_n}\rightarrow v^i$ in $C_{loc}^2(\R^{N+1})$ as ${\e_n}\rightarrow 0^+$.
Moreover, the elliptic $L^q$-estimate with $q>1$ yields $v^i\in C_{loc}^2(\R^{N+1})\cap W^{2,q}(\R^{N+1})$.
Obviously, $v^i\not\equiv0$ due to $u_\e(0)\geq c_1>0$.
Then using the similar argument as in Lemma \ref{Lem:upper}, $v^i$ solves equation (\ref{eqn:redu-qs}) with $d=|\bar{P}|$.
By applying the well-known Moser iteration, we can show that $|v^i|\in L^\infty(\R^N)$. Furthermore, by comparison principle, we obtain
\begin{equation}\label{eqn:zhu2--}
|v^i(x)|\leq C_1e^{-\frac{\sigma_k}{|\bar{P}|^{\al/2}}|x|}\quad\text{for}\,\,x\in\R^{N+1},\,\,i=1,...,k
\end{equation}
for some $\sigma_k>0$.
 For any given large constant $R$ in (\ref{eqn:zhu1}), we put
$\kappa_R:=\bar{C}_1e^{-\frac{\sigma_k R}{2}}$ for some $\bar{C}_1>0$. Then for $n$ large enough, we have
\begin{equation}\label{eqn:zhu2-}
\|v^i_{\e_n}-v^i\|_{C^2(\overline{B_{2R}(0)})}\leq \kappa_R.
\end{equation}
From (\ref{eqn:zhu1}), we deduce that
\begin{equation}\label{eqn:zhu2}
\aligned
c_{\e_n}&\geq\sum_{i=1}^{k}{\e_n}^{N+1}\int_{B_{R}(0)}\frac{1}{|{\e_n} x+P^i_{\e_n}|^\al}(\frac{1}{2}f(v^i_{\e_n})v^i_{\e_n}-F(v^i_{\e_n}))dx\\
&=\sum_{i=1}^{k}{\e_n}^{N+1}\bigg[\int_{B_{R}(0)}\frac{1}{|\bar{P}|^\al}(\frac{1}{2}f(v^i)v^i-F(v^i))dx+\mathcal{A}^i)\bigg]\\
&=\sum_{i=1}^{k}{\e_n}^{N+1}|\bar{P}|^{\frac{\al}{2}(N-1)}\bigg[\int_{B_{R|\bar{P}|^{-\al/2}}(0)}(\frac{1}{2}f(\bar{U}^i)\bar{U}^i-F(\bar{U}^i))dx+\mathcal{A}^i)\bigg]\\
&=\sum_{i=1}^{k}{\e_n}^{N+1}|\bar{P}|^{\frac{\al}{2}(N-1)}\bigg[\int_{\R^{N+1}}(\frac{1}{2}f(\bar{U}^i)\bar{U}^i-F(\bar{U}^i))dx
-\int_{D}(\frac{1}{2}f(\bar{U}^i)\bar{U}^i-F(\bar{U}^i))dx+\mathcal{A}^i)\bigg],
\endaligned
\end{equation}
where $D:=\R^{N+1}\setminus{B_{R\bar{P}^{-\al/2}}(0)}$, and
$\bar{U}^i(\cdot )=v^i(|\bar{P}|^{\frac{\al}{2}}\cdot)$ is a nontrivial solution of equation (\ref{eqn:spacefq}),
and
$$
\mathcal{A}^i=\int_{B_{R}(0)}\frac{1}{|{\e_n} x+P^i_{\e_n}|^\al}(\frac{1}{2}f(v^i_{\e_n})v^i_{\e_n}-F(v^i_{\e_n}))dx
-\int_{B_{R}(0)}\frac{1}{|\bar{P}|^\al}(\frac{1}{2}f(v^i)v^i-F(v^i))dx,
$$
which can be rewritten as
$$
\aligned
\mathcal{A}^i=&\int_{B_{R}(0)}\bigg(\frac{1}{|{\e_n} x+P^i_{\e_n}|^\al}-\frac{1}{|P^i_{\e_n}|^\al}\bigg)(\frac{1}{2}f(v^i_{\e_n})v^i_{\e_n}-F(v^i_{\e_n}))dx\\
&+\int_{B_{R}(0)}\bigg(\frac{1}{|P^i_{\e_n}|^\al}-\frac{1}{|\bar{P}|^\al}\bigg)(\frac{1}{2}f(v^i_{\e_n})v^i_{\e_n}-F(v^i_{\e_n}))dx\\
&+\int_{B_{R}(0)}\frac{1}{|\bar{P}|^\al}\bigg(\frac{1}{2}f(v^i_{\e_n})v^i_{\e_n}-F(v^i_{\e_n})-\frac{1}{2}f(v^i)v^i+F(v^i)\bigg)dx\\
:=& \mathcal{A}^i_1(\e_n)+\mathcal{A}^i_2(\e_n)+\mathcal{A}^i_3(\e_n).
\endaligned
$$
Since $v_{\e_n}$ is uniformly bounded and $f$ is of class $C^1$, by the mean-value theorem, one has
\begin{equation}\label{eqn:zhu2+}
|\mathcal{A}^i_1(\e_n)|\leq C\int_{B_{R}(0)}\bigg|\frac{1}{|{\e_n} x+P^i_{\e_n}|^\al}-\frac{1}{|P^i_{\e_n}|^\al}\bigg|dx\leq C_R\e_n.
\end{equation}
By the dominate convergence theorem, we also have
\begin{equation}\label{eqn:zhu2++}
|\mathcal{A}^i_2(\e_n)|\leq C\int_{B_{R}(0)}\bigg|\frac{1}{|P^i_{\e_n}|^\al}-\frac{1}{|\bar{P}|^\al}\bigg|dx\rightarrow0,\quad n\rightarrow\infty.
\end{equation}
From (\ref{eqn:zhu2-}) and the fact that $v_{\e_n}$ is uniformly bounded and $f$ is of class $C^1$, we deduce that for $n$ large enough,
\begin{equation}\label{eqn:zhu2+++}
|\mathcal{A}^i_3(\e_n)|\leq C|B_R(0)|\kappa_R=C|B_R(0)|\bar{C}_1e^{-\frac{\sigma_k R}{2}}.
\end{equation}
It follows from (\ref{eqn:zhu2--}) that there
exists $C_2>0$ such that
\begin{equation}\label{eqn:zhu3}
\int_{D}(\frac{1}{2}f(\bar{U}^i)\bar{U}^i-F(\bar{U}^i))dx\leq C_2e^{-\mu_k R}
\end{equation}
for some $\mu_k>0$ and $C_2>0$. Combining (\ref{eqn:zhu2})-(\ref{eqn:zhu3}) we obtain
\begin{equation}\label{eqn:zhu4}
\aligned
&\liminf\limits_{n\rightarrow\infty}{\e_n}^{-(N+1)}c_{\e_n}\\
\geq&\liminf\limits_{n\rightarrow\infty}
k|\bar{P}|^{\frac{\al}{2}(N-1)}\bigg[J(\bar{U})-C_R\e_n-|\mathcal{A}^i_3(\e_n)|-C|B_R(0)|\bar{C}_1e^{-\frac{\sigma_k R}{2}}
- C_2e^{-\mu_k R}\bigg]\\
=&
k|\bar{P}|^{\frac{\al}{2}(N-1)}\bigg[J(\bar{U})-C|B_R(0)|\bar{C}_1e^{-\frac{\sigma_k R}{2}}
- C_2e^{-\mu_k R}\bigg].
\endaligned
\end{equation}
Since $ C, \bar{C}_1, C_2>0$ in (\ref{eqn:zhu4}) are independent of $R$ and $\e_n$ for $n$ large enough,
and $U$ is the least energy positive solution,  from (\ref{eqn:zhu4}) we have
\begin{equation}\label{eqn:zhu5}
{\e_n}^{-(N+1)}c_{\e_n}\geq k
|\bar{P}|^{\frac{\al}{2}(N-1)}\bigg[J(U)-C|B_R(0)|\bar{C}_1e^{-\frac{\sigma_k R}{2}}
- C_2e^{-\mu_k R}\bigg].
\end{equation}
Choose
$$
k >\frac{2[d_1^{(N-1)\frac{\al}{2}}+d_2^{(N-1)\frac{\al}{2}}]}{|\bar{P}|^{\frac{\al}{2}(N-1)}}+1,
$$
 then letting $R\rightarrow\infty$, (\ref{eqn:zhu5})
implies a contradiction with Lemma \ref{Lem:upper}. Therefore, $P_\e\not\in \Omega^o$. Similarly, $Q_\e\not\in \Omega^o$.
We now show $\bar{P}\not\in\Omega^o$. Assume on the contrary that there exists $\{\e_n\}$ satisfying $\e_n\rightarrow 0^+$ such that
\begin{equation}\label{eqn:limit0}
P_{\e_n}=(\mathbf{x}_{\e_n},x_{\e_n,N+1})\in\Omega^1\rightarrow \bar{P}\in\Omega^o
\end{equation}
with $|\mathbf{x}_{\e_n}|=d_{\e_n}$ and such that $d_{\e_n}\rightarrow d$ for some $d\geq\frac{a^2}{2}$,
 as $n\rightarrow +\infty$. We can use the almost same arguments as above to get a contradiction. Therefore,
 $\bar{P}\in\Omega^1$ and $\bar{Q}\in\Omega^1$, which imply by the symmetric of $u_\e\in H_{\sharp}(\Omega)$ that
$$
 \frac{\min\{dist(\pm P_\e,\Omega^o), dist(\pm Q_\e,\Omega^o)\}}{\e}\rightarrow\infty,
 $$
 as $\e\rightarrow\infty$. The proof is complete.
\qed

Define function space
$$
\aligned
&H^+:=\{u\in H_{\sharp}(\Omega)\cap H^1(\Omega^1):\,\,\frac{\partial u}{\partial x_{N+1}}\bigg|_{\Omega^o}=0\},\\
&H^-:=\{u\in H_{\sharp}(\Omega)\cap H^1(\Omega^2):\,\,\frac{\partial u}{\partial x_{N+1}}\bigg|_{\Omega^o}=0\},
\endaligned
$$
which will be used in the following lemma.

Without loss of generality,
we assume $P_\e,Q_\e\in \Omega^1$, then we have the following results.

\begin{lemma}\label{Lem:asym} As $\e\rightarrow0$, we have
$\frac{|P_\e-Q_\e|}{\e}\rightarrow\infty$ as $\e\rightarrow0$.
\end{lemma}
\Proof
Assume on the contrary that $\frac{|P_\e-Q_\e|}{\e}\rightarrow\kappa<\infty$. Let us show first that $\kappa>0$.
Since ${\e}^2\triangle u_{\e}=\frac{1}{|x|^\al}(u_{\e}-f(u_{\e}))$ and ${\e}^2\triangle u_{\e}(P_{\e})\leq0$, we have
$u_{\e}(P_{\e})\leq f(u_{\e}(P_{\e}))$ which implies by condition ($f_1$) that there exists $c_1>0$ such that
$u_\e(P_\e)\geq c_1$. Similarly, we also have $u_\e(Q_\e)\leq -c_2$ for some $c_2>0$. Thus,
we have $u_\e(P_\e)-u_\e(Q_\e)\geq c_1+c_2$.  If there exists one point $R_\e\in \overline{P_\e Q_\e} \cap \partial \Omega^1$,
then by recalling Lemma \ref{Lem:axis}, one has
$$
\min\{|R_\e-P_\e|,|R_\e-Q_\e|\}/\e\rightarrow\infty,\quad\quad\text{as}\,\,\e\rightarrow0^+,
$$
which implies the conclusion of Lemma \ref{Lem:asym} holds true.
Thus, segment $\overline{P_\e Q_\e}\subset \Omega^1$.
From the boundedness of $\{u_\e\}$
in $L^\infty(\Omega)$ and Schauder's estimates, we deduce that $\e|\nabla u_\e|\leq C$ for some $C>0$
independently of $\varepsilon$. Thus,
$$
c_1+c_2\leq |u_\e(P_\e)-u_\e(Q_\e)|\leq \e|\nabla u_\e(\xi)|\frac{|P_\e-Q_\e|}{\e}
$$
for some $\xi\in \Omega$.
 Therefore, $\kappa>0$ and $O:=\lim_{\e\rightarrow0}\frac{Q_\e-P_\e}{\e}\not\in\R^{N+1}\setminus\{0\}$.\\
Let us define $u_\e=\bar{u}_\e+\tilde{u}_\e$ by
$$
u_\e=\left\{
  \begin{array}{ll}
  \displaystyle  \bar{u}_\e&  \,\, x\in\Omega^1,\\
   \tilde{u}_\e,& \,\, x\in\Omega^2,
  \end{array}
\right.
$$
then since $u_\e\in H_{\sharp}(\Omega)$, $\bar{u}_\e$ satisfies the following equation
\begin{equation}\label{eqn:half-eq1}
\left\{
  \begin{array}{ll}
  \displaystyle -\e^2\triangle u+\frac{u}{|x|^\al}=\frac{f(u)}{|x|^\al}&  \,\, x\in\Omega^1,\\
   \frac{\partial u}{\partial n}=0,& \,\, x\in\Omega^o,\\
    u=0,&\,\,x\in\partial\Omega^1\setminus{\Omega^o},
  \end{array}
\right.
\end{equation}
whose energy functional $J_{\e,\Omega^1}: H^+:\rightarrow\R$ defined by
$$
J_{\e,\Omega^1}(u)=\frac{\varepsilon^2}{2}\int_{\Omega^1}|\nabla u|^2dx+\frac{1}{2}\int_{\Omega^1}\frac{u^2}{|x|^\alpha}dx
-\int_{\Omega^1}\frac{F(u)}{|x|^{\alpha}}dx.
$$
We have similarly that $\tilde{u}_\e$ satisfies
\begin{equation}\label{eqn:half-eq2}
\left\{
  \begin{array}{ll}
  \displaystyle -\e^2\triangle u+\frac{u}{|x|^\al}=\frac{f(u)}{|x|^\al}&  \,\, x\in\Omega^2,\\
   \frac{\partial u}{\partial n}=0,& \,\, x\in\Omega^o,\\
    u=0,&\,\,x\in\partial\Omega^2\setminus{\Omega^o},
  \end{array}
\right.
\end{equation}
whose energy functional $J_{\e,\Omega^2}: H^-:\rightarrow\R$ defined by
$$
J_{\e,\Omega^2}(u)=\frac{\varepsilon^2}{2}\int_{\Omega^2}|\nabla u|^2dx+\frac{1}{2}\int_{\Omega^2}\frac{u^2}{|x|^\alpha}dx
-\int_{\Omega^2}\frac{F(u)}{|x|^{\alpha}}dx.
$$
By the definition of $H_{\sharp}(\Omega)$ and the fact that functional $J_\epsilon(u)$ is even at $u$, we immediately get
$J_{\e,\Omega^1}(\bar{u}_\e)=J_{\e,\Omega^2}(\tilde{u}_\e)$ and $J_\e(u_\e)=J_{\e,\Omega^1}(\bar{u}_\e)+J_{\e,\Omega^2}(\tilde{u}_\e)$.
It follows from (\ref{eqn:nahar7}) that
\begin{equation}\label{eqn:half-3}
J_{\e}(\bar{u}_\e)\leq \e^{(N+1)}(d_1^{(N-1)\frac{\al}{2}}+d_2^{(N-1)\frac{\al}{2}})[J(U)+O(\e)].
\end{equation}
Set $w_\e(y)=\bar{u}_\e(\e y+P_\e)$, then by the elliptic regularity theory, we can easily show that
$w_\e(y)\rightarrow w$ in $C_{loc}^2(\R^{N+1})$ and $w$ satisfies
\begin{equation}\label{eqn:half-4}
\left\{
  \begin{array}{ll}
  \displaystyle  -\triangle {w}+\frac{w}{|\bar{P}|^\alpha} =\frac{f(w)}{|\bar{P}|^\alpha}&  \mbox{in}\,\, \R^{N+1},\\
   w(0)\geq c_1,\,\, w(O)\leq-c_2,
  \end{array}
\right.
\end{equation}
which implies $w$ is a nodal solution of equation (\ref{eqn:half-4}) whose energy functional is
$$
J_{\bar{P}}(w)=\frac{1}{2}\int_{\R^{N+1}}|\nabla w|^2+\frac{1}{2}\int_{\R^{N+1}}\frac{w^2}{|\bar{P}|^\al}dx
-\int_{\R^{N+1}}\frac{F(w)}{|\bar{P}|^\al}dx.
$$
Thus, we have
$$
\int_{\R^{N+1}}|\nabla w^\pm|^2+\int_{\R^{N+1}}\frac{|w^\pm|^2}{|\bar{P}|^\al}dx
=\int_{\R^{N+1}}\frac{f(w^\pm)w^\pm}{|\bar{P}|^\al}dx.
$$
From the above facts, we deduce that
\begin{equation}\label{eqn:half-4+}
\aligned
J_{\bar{P}}(w^\pm)&=\int_{\R^{N+1}}\frac{1}{|\bar{P}|^\al}[\frac{1}{2}f(w^\pm)w^\pm-F(w^\pm)]dx\\
&\geq \int_{\R^{N+1}}\frac{1}{|\bar{P}|^\al}[\frac{1}{2}f(W)W-F(W)]dx=J_{\bar{P}}(W),
\endaligned
\end{equation}
where $W$ is the unique positive solution of (\ref{eqn:redu-qs}) with $d=|\bar{P}|$. Without loss of generality, assume
$|\bar{Q}|\leq |\bar{P}|$, then by (\ref{eqn:half-3}) we have
$$
J_{\e,\Omega^1}(\bar{u}_\e)\leq 2\e^{N+1}|\bar{P}|^{(N-1)\frac{\al}{2}}[J(U)+O(\e)],
$$
which implies by Fatou's lemma that
\begin{equation}\label{eqn:half-5}
\aligned
J_{\bar{P}}(w)&=\int_{\R^{N+1}}\frac{1}{|\bar{P}|^\al}[\frac{1}{2}f(w)w-F(w)]dx\\
&\leq \liminf\limits_{\e\rightarrow0}\e^{-(N+1)}[J_{\e,\Omega^1}(\bar{u}_\e)-\frac{1}{2}J'_{\e,\Omega^1}(\bar{u}_\e)\bar{u}_\e]\\
&\leq 2\int_{\R^{N+1}}\frac{1}{|\bar{P}|^\al}[\frac{1}{2}f(W)W-F(W)]dx\\
&=2J_{\bar{P}}(W).
\endaligned
\end{equation}
Combining (\ref{eqn:half-4+}) with (\ref{eqn:half-5}), we have
$J_{\bar{P}}(W)=J_{\bar{P}}(w^\pm)$. Since $W$ is the unique positive solution of (\ref{eqn:redu-qs}) with $d=|\bar{P}|$,
$W$ minimizes the functional $J_{\bar{P}}$ on the following Nehari manifold
$$
\mathcal{N}_{\bar{P}}=\{u\in H_0^1(\Omega):\,\,u\not\equiv 0,\,J'_{\bar{P}}(u)u=0\}.
$$
Clearly, $w^\pm\in \mathcal{N}_{\bar{P}}$ and minimizes $J_{\bar{P}}$ on $\mathcal{N}_{\bar{P}}$. Hence, it is easy to show that
$w^\pm$ are also positive solutions of (\ref{eqn:redu-qs}) with $d=|\bar{P}|$. By the uniqueness of positive solution, we obtain
$w^+(x)=W(x-\bar{x}_1)$ and $w^-(x)=-W(x-\bar{x}_2)$ for some $\bar{x}_1,\bar{x}_2\in\R^{N+1}$.
Therefore, $w^-(x)<0$ for any $x\in \R^{N+1}$, which is a contradiction. The proof is complete.
\qed

In Lemma \ref{Lem:upper}, we replace $P$ and $Q$ by $P_\e$ and $Q_\e$, respectively.
Based on Lemma \ref{Lem:axis} and Lemma \ref{Lem:asym}, similar arguments as Lemma \ref{Lem:upper} yield
\begin{equation}\label{eqn:ex-upper}
c_\e \leq2\e^{N+1}(|\bar{P}|^{(N-1)\frac{\al}{2}}+|\bar{Q}|^{(N-1)\frac{\al}{2}})[J(U)+O(\e)].
\end{equation}
Assume $u_\e\in H_\sharp(\Omega)$ is the nodal solution obtained in Theorem \ref{Thm:exist}. We now show
\begin{lemma}\label{Lem:limit}
$\lim\limits_{\e\rightarrow0^+}\e^{-(N+1)}J_\e(u_\e)=2(\bar{P}^{(N-1)\frac{\al}{2}}+\bar{Q}^{(N-1)\frac{\al}{2}})J(U)$.
\end{lemma}
\Proof
 Let us define
$$
r_\e=\min\bigg\{dist(P_\e,\partial\Omega),dist(Q_\e,\partial\Omega)),\frac{|P_\e-Q_\e|}{2}\bigg\},
$$
Then by recalling Lemma \ref{Lem:axis}, Lemma \ref{Lem:asym}, and the definition of $H_{\sharp}(\Omega)$, we can deduce that $r_\e/\e\rightarrow+\infty$
as $\e\rightarrow0^+$,
and
$$
B_{r_\e}(P_\e), B_{r_\e}(-P_\e), B_{r_\e}(Q_\e), B_{r_\e}(-Q_\e)\subset \Omega,
$$
which are disjoint each other. Set
$$
\mathcal{D}:=B_{r_\e}(P_\e)\cup B_{r_\e}(-P_\e)\cup B_{r_\e}(Q_\e)\cup B_{r_\e}(-Q_\e).
$$
Since $u_\e\in H_{\sharp}(\Omega)$ is a sign-changing solution of equation (\ref{eqn:redu-qs-}), by ($f_2$) we have
\begin{equation}\label{eqn:guji1}
\aligned
&\e^{-(N+1)}J_{\varepsilon}(u_{\varepsilon})\\
&\geq\e^{-(N+1)}\int_{\mathcal{D}}\frac{1}{|x|^\al}(\frac{1}{2}f(u_{\e})u_{\e}-F(u_{\e}))dx\\
&=\e^{-(N+1)}\int_{B_{r_\e}(P_\e)}\frac{1}{|x|^\al}(\frac{1}{2}f(u_{\e})u_{\e}-F(u_{\e}))dx+\e^{-(N+1)}\int_{B_{r_\e}(-P_\e)}\frac{1}{|x|^\al}(\frac{1}{2}f(u_{\e_n})u_{\e}-F(u_{\e}))dx\\
&+\e^{-(N+1)}\int_{B_{r_\e}(Q_\e)}\frac{1}{|x|^\al}(\frac{1}{2}f(u_{\e})u_{\e}-F(u_{\e}))dx+\e^{-(N+1)}\int_{B_{r_\e}(-Q_\e)}\frac{1}{|x|^\al}(\frac{1}{2}f(u_{\e_n})u_{\e}-F(u_{\e}))dx\\
&=\mathcal{B}_1+\mathcal{B}_2+\mathcal{B}_3+\mathcal{B}_4.
\endaligned
\end{equation}
By making the change of variable $v_\e(x)=u_\e(\e x+P_\e)$ for $x\in B_{r_\e}(P_\e)$,
using elliptic regularity theory and the boundedness of $\{u_\e\}$ in $H_{\sharp}(\Omega)$,
we obtain
$v_\e\rightarrow v$ in $C_{loc}^2(\R^{N+1})$ as $\e\rightarrow 0^+$. Moreover, the elliptic $L^q$-estimate with $q>1$ yields $v\in C_{loc}^2(\R^{N+1})\cap W^{2,q}(\R^{N+1})$.
We note that $v$ is a nontrivial solution of equation
(\ref{eqn:redu-qs}) with $d=|\bar{P}|$, and $\bar{U}(\cdot)=v(|\bar{P}|^{\al/2}\cdot)$ is a nontrivial solution of equation (\ref{eqn:spacefq}).
Based on the above facts, by the Fatou's lemma, we get
\begin{equation}\label{eqn:guji2}
\aligned
\liminf\limits_{\e\rightarrow0^+}\mathcal{B}_1&= \liminf\limits_{\e\rightarrow0^+}\int_{B_{\frac{r_\e}{\e}}(P_\e)}\frac{1}{|\e x+P_{\e}|^\al}(\frac{1}{2}f(v_{\e})v_{\e}-F(v_{\e}))dx\\
&\geq\int_{\R^{N+1}}\frac{1}{|\bar{P}|^\al}(\frac{1}{2}f(v)v-F(v))dx\\
&\geq|\bar{P}|^{\frac{\al}{2}(N-1)}J(\bar{U})\\
&\geq|\bar{P}|^{\frac{\al}{2}(N-1)}J(U),
\endaligned
\end{equation}
where we use the fact that $J(\bar{U})\geq J(U)$ since $U\in H^1(\R^{N+1})$ is the unique positive solution of (\ref{eqn:spacefq}).
Arguing similarly as above, we also have
\begin{equation}\label{eqn:guji3}
\aligned
&\liminf\limits_{\e\rightarrow0^+}\mathcal{B}_2\geq|\bar{P}|^{\frac{\al}{2}(N-1)}J(U)\\
&\liminf\limits_{\e\rightarrow0^+}\mathcal{B}_3\geq |\bar{Q}|^{\frac{\al}{2}(N-1)}J(U)\\
&\liminf\limits_{\e\rightarrow0^+}\mathcal{B}_4\geq|\bar{Q}|^{\frac{\al}{2}(N-1)}J(U) .
\endaligned
\end{equation}
Since functional $J$ is even, by combining (\ref{eqn:guji1})-(\ref{eqn:guji3}), we get immediately
\begin{equation}\label{eqn:guji5}
\liminf\limits_{\e\rightarrow0^+}\e^{-(N+1)}J_\e(u_\e)
\geq2(\bar{P}^{(N-1)\frac{\al}{2}}+\bar{Q}^{(N-1)\frac{\al}{2}})J(U).
\end{equation}
 Therefore, the conclusion follows immediately from (\ref{eqn:ex-upper}) and (\ref{eqn:guji5}).
 The proof is complete.
\qed

Recall that $\bar{u}_\e(x)$ has been defined in Lemma \ref{Lem:asym},
$\bar{u}_\e(x)=u_\e(x)$ for $x\in \Omega^1$ and satisfies equation (\ref{eqn:half-eq1}).
Note that $J_{\e,\Omega^1}(\bar{u}_\e)=\frac{1}{2}c_\e$, where $c_\e$ has been defined in Theorem \ref{Thm:exist}.
Let us define the least sign changing energy level of $J_{\e,\Omega^1}$  as follows
$$
c_\e^+:=\inf\limits_{u\in \mathcal{N}_\e^+} J_{\e,\Omega^1}(u),\quad\quad  \mathcal{N}_\e^+:=\{u\in H^+:\,\,u^{\pm}\not=0,\,J'_{\e,\Omega^1}(u^+)u^+=J'_{\e,\Omega^1}(u^-)u^-=0\}.
$$
\begin{lemma}\label{Lem:half-enery}
$c_\e^+=\frac{1}{2}c_\e$.
\end{lemma}
\Proof
Using the similar arguments as in Theorem \ref{Thm:exist}, we can prove that there exists a least energy sign changing solution
$u\in H^+$ of equation (\ref{eqn:half-eq1}) with $J_{\e,\Omega^1}(u)=c_\e^+$.
Since $\bar{u}_\e$ is a nodal solution of equation (\ref{eqn:half-eq1}), then $\frac{c_\e}{2}\geq c_\e^+$.
Define $v=u(\mathbf{x},-x_{N+1})$ for $x\in\Omega^2$, then $v$ is a sign-changing solution of equation (\ref{eqn:half-eq2}),
and then by symmetry, $J_{\e,\Omega^2}(v)=c_\e^+$. Define
$$
\mathfrak{u}(x):=\left\{
  \begin{array}{ll}
  \displaystyle  u(x)&  \,\, x\in\Omega^1,\\
   v(x),& \,\, x\in\Omega^2,
  \end{array}
\right.
$$
then $\mathfrak{u}\in H_{\sharp}(\Omega)$ and
$$
J'_{\e}(\mathfrak{u})=0,\quad \text{and}\quad J_{\e}(\mathfrak{u})=J_{\e,\Omega^2}(v)+J_{\e,\Omega^1}(u)=2c_\e^+.
$$
Based on the definition of $c_\e$, we have $c_\e^+\geq \frac{c_\e}{2}$. The proof is complete.
\qed

Let us define
$$
\Omega^1_{P_\e}:=\{y\in\R^{N+1}|\,\,\e y+P_\e\in \Omega_{\e,+}^1\},
\quad\Omega^1_{Q_\e}:=\{y\in\R^{N+1}|\,\,\e y+ Q_\e\in \Omega_{\e,-}^1\},
$$
where $\Omega_{\e,\pm}^1$ are the support sets of $\bar{u}_\e^\pm$, respectively.

\begin{lemma}\label{Lem:limit-blow}
As $\e\rightarrow0$, we have $\Omega^1_{P_\e}\rightarrow\R^{N+1}$ and $\Omega^1_{Q_\e}\rightarrow\R^{N+1}$.
\end{lemma}

\Proof
We first show that $\Omega_{\e,\pm}^1$ are both connected domains. Assume on the contrary that the number of nodal domains of
$\bar{u}_\e$ is bigger than $2$. Let us denote nodal domain by $\{\Omega^1_i\}, i=1,2,...,k$ with $k\geq 3$.
Let us define
$\bar{u}_{\e}=\bar{u}_{\e,1}+\bar{u}_{\e,2}+\bar{u}_{\e,3}$
with
$$
\bar{u}_{\e,i}\not=0,\,\, \bar{u}_{\e,1}\geq0,\,\,\bar{u}_{\e,2}\leq0\quad \text{and}\quad suppt(\bar{u}_{\e,i})
\cap suppt(\bar{u}_{\e,j})=\emptyset,
\quad\text{for}\,\,i\not=j,\,\,i,j=1,2,3.
$$
Then it is easy to see that
$$
J'_{\e,\Omega^1}(\bar{u}_{\e,i})\bar{u}_{\e,i}=0,\quad \text {for}\,\,i=1,2,3.
$$
By ($f_2$), we have $J_{\e,\Omega^1}(\bar{u}_{\e,i})>0$ for all
$i=1,2,3$.
In virtue of Lemma \ref{Lem:half-enery} and the definition of $\mathcal{N}_\e^+$, we have
$$
\frac{1}{2}c_\e\leq J_{\e,\Omega^1}(\bar{u}_{\e,1}+\bar{u}_{\e,2})<\sum\limits_{i=1}^3J_{\e,\Omega^1}(\bar{u}_{\e,i})
=J_{\e,\Omega^1}(\sum\limits_{i=1}^3\bar{u}_{\e,i})=\frac{1}{2}c_\e.
$$
It is a contradiction. It is easy to see from the classical elliptic regularity theory that ${\bar{u}}_\e\in C^2(\Omega^1)$,
which implies that $\Omega_{\e,\pm}^1$ are both connected domains.

Without loss of generality, we assume $\Omega^1_{P_\e}\rightarrow \Omega^1_{\bar{P}}$ and $\Omega^1_{Q_\e}\rightarrow \Omega^1_{\bar{Q}}$
as $\e\rightarrow0^+$.
Now we only need to prove $\Omega^1_{\bar{P}}=\Omega^1_{\bar{Q}}=\R^{N+1}$.
Let us define the part boundary of nodal domain $\Omega^1_{\e,+}$ as follows
$$
\Theta:=\{x\in\partial\Omega_{\e,+}^1,\,\, x\not\in \Omega^o\,\,\text{and}\,\,x\not\in \partial\Omega\}.
$$
Since  $\Omega_{\e,\pm}^1$ are both connected domains, $\Theta$ is also the part boundary of nodal domain $\Omega^1_{\e,-}$.
Obviously, $\bar{u}^+_\e|_\Theta=\bar{u}^-_\e|_\Theta=0$.
We first show
\begin{equation}\label{eqn:blow3}
\frac{\min\{dist(P_\e,\Theta), dist(Q_\e,\Theta)\}}{\e}\rightarrow\infty.
\end{equation}
as $\e\rightarrow0^+$.
We only prove that $dist(P_\e,\Theta)/{\e}\rightarrow\infty$, the remaining part can be obtained similarly.
Assume on the contrary that there exists $\kappa_1\geq0$ such that $\lim_{\e\rightarrow0}dist(P_\e,\Theta)/\e=\kappa_1$.
Thus, we may assume that $P_\e\rightarrow \bar{P}\in \Theta$ and there exists $\bar{P}_\e\in \Theta$ such that
$$
dist(P_\e,\Theta)=|P_\e-\bar{P}_\e|,\quad\text{and}\quad \bar{P}_\e\rightarrow \bar{P}\,\,\text{as}\,\,\e\rightarrow0^+.
$$
Then we have $\frac{|P_\e-\bar{P}_\e|}{\e}\rightarrow\kappa_1<\infty$.
Since it is easy to obtain from the Maximum principle that there exists $c_1>0$ such that
$u_\e(P_\e)\geq c_1$. Then by $u_\e(\bar{P}_\e)=0$, we have $u_\e(P_\e)-u_\e(\bar{P}_\e)\geq c_1$.
Arguing similarly as in Lemma \ref{Lem:asym}, we have $\kappa_1>0$.
 Set $\bar{O}:=\lim_{\e\rightarrow0}\frac{\bar{P}_\e-P_\e}{\e}\in\R^{N+1}\setminus\{0\}$.
Set $w_\e(y)=\bar{u}_\e(\e y+P_\e)$, then by Lemma \ref{Lem:axis} and the elliptic regularity theory, we can easily show that
$w_\e(y)\rightarrow w$ in $C_{loc}^2(\R^{N+1})$ as $\e\rightarrow0$, and $w$ satisfies
\begin{equation}\label{eqn:blow4}
\left\{
  \begin{array}{ll}
  \displaystyle  -\triangle {w}+\frac{w}{|\bar{P}|^\alpha} =\frac{f(w)}{|\bar{P}|^\alpha}&  \mbox{in}\,\, \R^{N+1},\\
   w(0)\geq c_1,\,\, w(\bar{O})=0.
  \end{array}
\right.
\end{equation}
By the well-known strong Maximum principle, we deduce that $w$ is a nodal solution of equation (\ref{eqn:blow4}).
Hence, using the similar arguments as Lemma \ref{Lem:asym}, we have $J_{\bar{P}}(W_{\bar{P}})=J_{\bar{P}}(w^\pm)$,
where $W_{\bar{P}}$ is the unique positive solution of (\ref{eqn:redu-qs}) with $d=|\bar{P}|$. It implies a contradiction.
Then, (\ref{eqn:blow3}) holds true.
Using the fact that $\Omega_{\e,\pm}^1$ are both connected domains, combining (\ref{eqn:blow3}) and Lemma \ref{Lem:axis}, we get immediately
$$
\frac{\min\{dist(P_\e,\partial\Omega^1_{\e,+}), dist(Q_\e,\partial\Omega^1_{\e,-})\}}{\e}\rightarrow\infty,
\quad\text{as}\,\,\,\e\rightarrow0^+,
$$
from which we deduce immediately the conclusions of lemma. The proof is complete.
\qed

\begin{lemma}\label{Lem:Unique}
Assume $u_\e\in H_\sharp(\Omega)$ is the solution obtained in Theorem \ref{Thm:exist}. Then $u_\e$ has only two
positive local maximums and only two negative local minimums.
\end{lemma}
\Proof
Since $u_\e\in H_\sharp(\Omega)$ and $P_\e$ is a local maximum point,
the antipodal point $-P_\e$ is also a maximum point of $u_\e$. Correspondingly,
$Q_\e$ and $-Q_\e$ are minimum points of $u_\e$. Assume $\bar{u}_\e$ solves equation (\ref{eqn:half-eq1}) which has been defined in Lemma \ref{Lem:asym}.
Let us define
$\bar{v}_\e^+(y):=\bar{u}_\e^+(\e y+P_\e)$ and $\bar{v}_\e^-(y):=\bar{u}_\e^-(\e y+Q_\e)$,
then by Remark 1.1, we have
\begin{equation}\label{eqn:decay1}
\aligned
&\e^{-(N+1)}J_{\e,\Omega^1}(\bar{u}_\e)\\
=&\e^{-(N+1)}J_{\e,\Omega^1}(\bar{u}^+_\e)+\e^{-(N+1)}J_{\e,\Omega^1}(\bar{u}^-_\e)\\
=&\e^{-(N+1)}\frac{\mu-1}{2\mu+2}\bigg[\int_{\Omega_{\e,+}^1}(|\nabla \bar{u}^+_\e|^2+\frac{|\bar{u}^+_\e|^2}{|x|^\al})dx
+\int_{\Omega_{\e,+}^1}(|\nabla \bar{u}^-_\e|^2+\frac{|\bar{u}^-_\e|^2}{|x|^\al})dx\bigg]\\
=&\frac{\mu-1}{2\mu+2}\bigg[\int_{\Omega^1_{P_\e}}(|\nabla \bar{v}^+_\e|^2+\frac{|\bar{v}^+_\e|^2}{|\e x+P_\e|^\al})dx
+\int_{\Omega^1_{Q_\e}}(|\nabla \bar{v}^-_\e|^2+\frac{|\bar{v}^-_\e|^2}{|\e x+Q_\e|^\al})dx\bigg]\\
=&J_{\e,\Omega^1_{P_\e}}(\bar{v}^+_\e)+J_{\e,\Omega^1_{Q_\e}}(\bar{v}^-_\e).
\endaligned
\end{equation}
Here, $J_{\e,\Omega^1_{P_\e}}$ and $J_{\e,\Omega^1_{Q_\e}}$,  denote
functional $J_{\e}$ whose integral region are $\Omega^1_{P_\e}$ and $\Omega^1_{Q_\e}$, respectively.
It is easy to see from Lemma \ref{Lem:limit} and (\ref{eqn:decay1}) that $\|\bar{v}^+_\e\|_{H^1(\Omega^1_{P_\e})}$ and
$\|\bar{v}^-_\e\|_{H^1(\Omega^1_{P_\e})}$ are bounded uniformly for $\e$. Using the elliptic regularity theory and Lemma \ref{Lem:limit-blow}, we deduce that
\begin{equation}\label{eqn:decay2}
\bar{v}^+_\e\rightarrow \bar{v}^+\quad\text{in}\,\, C^2_{loc}(\R^{N+1})\quad\quad\text{and}\quad\quad
\bar{v}^-_\e\rightarrow \bar{v}^-\quad\text{in}\,\, C^2_{loc}(\R^{N+1})
\end{equation}
as $\e\rightarrow0^+$. Obviously, $\bar{v}^+$ and $\bar{v}^-$ are nontrivial solutions of equation (\ref{eqn:redu-qs})
with $d=|\bar{P}|$ and $d=|\bar{Q}|$, respectively.
In virtue of Lemma \ref{Lem:limit}, \ref{Lem:half-enery} and \ref{Lem:limit-blow}, from (\ref{eqn:decay1}), Fatou's lemma and ($f_2$), we deduce that
\begin{equation}\label{eqn:decay3}
\aligned
&(\bar{P}^{(N-1)\frac{\al}{2}}+\bar{Q}^{(N-1)\frac{\al}{2}})J(U)\\
&=\liminf\limits_{\e\rightarrow0}\bigg[\e^{-(N+1)}J_{\e,\Omega^1}(\bar{u}^+_\e)+\e^{-(N+1)}J_{\e,\Omega^1}(\bar{u}^-_\e)\bigg]\\
&=\liminf\limits_{\e\rightarrow0}\frac{\mu-1}{2\mu+2}\bigg[\int_{\Omega^1_{P_\e}}(|\nabla \bar{v}^+_\e|^2+\frac{|\bar{v}^+_\e|^2}{|\e x+P_\e|^\al})dx
+\int_{\Omega^1_{Q_\e}}(|\nabla \bar{v}^-_\e|^2+\frac{|\bar{v}^-_\e|^2}{|\e x+Q_\e|^\al})dx\bigg]\\
&\geq\frac{\mu-1}{2\mu+2}\bigg[\int_{\R^{N+1}}(|\nabla \bar{v}^+|^2+\frac{|\bar{v}^+|^2}{|\bar{P}|^\al})dx
+\int_{\R^{N+1}}(|\nabla \bar{v}^-|^2+\frac{|\bar{v}^-|^2}{|\bar{Q}|^\al})dx\bigg]\\
&=J_{\bar{P}}(\bar{v}^+)+J_{\bar{Q}}(\bar{v}^-)\\
&\geq J_{\bar{P}}(W_{\bar{P}})+J_{\bar{Q}}(W_{\bar{Q}})\\
&=(\bar{P}^{(N-1)\frac{\al}{2}}+\bar{Q}^{(N-1)\frac{\al}{2}})J(U).
\endaligned
\end{equation}
Here $J_{\bar{P}}$ and $J_{\bar{Q}}$ are the energy functional of equation (\ref{eqn:redu-qs})
with $d=|\bar{P}|$ and $d=|\bar{Q}|$, respectively.
So, $W_{\bar{P}}\in H^1(\R^{N+1})$ and $W_{\bar{Q}}\in H^1(\R^{N+1})$ are the least energy solutions
of equation (\ref{eqn:redu-qs}) with $d=|\bar{P}|$ and $d=|\bar{Q}|$, respectively.
By (\ref{eqn:decay3}) we have immediately
\begin{equation}\label{eqn:decay2-H}
\bar{v}^+_\e\rightarrow \bar{v}^+\quad\text{in}\,\,H^1(\R^{N+1})\quad\quad\text{and}\quad\quad
\bar{v}^-_\e\rightarrow \bar{v}^-\quad\text{in}\,\, H^1(\R^{N+1})
\end{equation}
as $\e\rightarrow0^+$.

We now prove that $\bar{u}_\e$ has at most one local maximum point. Assume on the contrary that
$\bar{u}_\e$ has two local maxima at $P_\e$ and $P'_\e$. Then there are three cases to be considered as follows:\\
Case 1: Suppose $\lim_{\e\rightarrow0}|P_\e-P'_\e|/\epsilon=0$. Observe by two local maximum points that
\begin{equation}\label{eqn:b1}
\nabla \bar{v}^+_\e(0)=\nabla \bar{v}^+_\e((P'_\e-P_\e)/\epsilon)=0.
\end{equation}
Since $\bar{v}^+$ is a positive solution of equation (\ref{eqn:redu-qs}) with $d=|\bar{P}|$, $\bar{v}^+(0)=\max_{\R^{N+1}}\bar{v}^+$
and $\bar{v}^+(|x|)=\bar{v}^+(r)$ is strictly decreasing at $r$, it is easy to obtain that $\triangle \bar{v}^+(0)<0$.
Therefore, by (\ref{eqn:decay2}), we also have $\triangle \bar{v}^+_\e(0)<0$ for small $\e$
and then
\begin{equation}\label{eqn:b2}
\triangle \bar{v}^+_\e(x)<0\quad \text{for}\,\, |x|\leq \varrho
\end{equation}
with $\varrho$ small enough.
In virtue of (\ref{eqn:b1}) and (\ref{eqn:b2}), it is easy to deduce that $P'_\e=P_\e$.\\
Case 2: Suppose $\lim_{\e\rightarrow0}|P_\e-P'_\e|/\epsilon=\beta>0$.
Assume $\tilde{O}:=\lim_{\e\rightarrow0}\frac{P'_\e-P_\e}{\e}\in \R^{N+1}\setminus\{0\}$, then
$\bar{v}^+>0$ satisfies
\begin{equation}\label{eqn:b3}
\left\{
  \begin{array}{ll}
  \displaystyle  -\triangle {\bar{v}^+}+\frac{\bar{v}^+}{|\bar{P}|^\alpha} =\frac{f(\bar{v}^+)}{|\bar{P}|^\alpha}&  \mbox{in}\,\, \R^{N+1},\\
   \nabla \bar{v}^+(0)=0,\, \nabla \bar{v}^+(\tilde{O})=0,
  \end{array}
\right.
\end{equation}
which contradicts with the fact that $\bar{v}^+(r)$ is strictly decreasing at $r$.\\
Case 3: Suppose $|P_\e-P'_\e|/\e\rightarrow\infty$ as $\e\rightarrow\infty$.
Assume $\lim_{\e\rightarrow0}P'_\e=\bar{P}'$.
Recalling Lemma \ref{Lem:asym} and Lemma \ref{Lem:limit},
we can get respectively $|Q_\e-P'_\e|/\e\rightarrow\infty$ and
$$
\frac{dist(P'_\e,\partial\Omega^1_{\e,+})}{\e}\rightarrow\infty,
\quad\text{as}\,\,\,\e\rightarrow0^+.
$$
Let us define
$$
r_\e=\min\bigg\{dist(P_\e,\partial\Omega^1_{\e,+}),dist(P'_\e,\partial\Omega^1_{\e,+})),dist(Q_\e,\partial\Omega^1_{\e,-})),\frac{|P_\e-P'_\e|}{2}\bigg\},
$$
then by Lemma \ref{Lem:limit}, we have
$r_\e/\e\rightarrow\infty$ as $\e\rightarrow0^+$. Set
$$
\bar{v}_{\e,1}(y):=\bar{u}_\e(\e y+P_\e),\quad \bar{v}_{\e,2}(y):=\bar{u}_\e(\e y+P'_\e),\quad
\bar{v}_{\e,3}(y):=\bar{u}_\e(\e y+Q_\e),\quad y\in B_{r_\e}(0).
$$
Using the elliptic regularity theory, we have
$$
\bar{v}_{\e,1}\rightarrow \bar{v}_{1}\quad\text{in}\,\, C^2_{loc}(\R^{N+1}),\quad\quad \bar{v}_{\e,2}\rightarrow \bar{v}_{2}\quad\text{in}\,\, C^2_{loc}(\R^{N+1}),
\quad\quad \bar{v}_{\e,3}\rightarrow \bar{v}_{3}\quad\text{in}\,\, C^2_{loc}(\R^{N+1}),
$$
as $\e\rightarrow0^+$. $\bar{v}_{1},\bar{v}_{2}$, and $\bar{v}_{3}$ are the positive solutions of equation
(\ref{eqn:redu-qs}) with $d=|\bar{P}|$, $d=|\bar{P'}|$, and $d=|\bar{Q}|$, respectively.
Using change of variables and Fatou's lemma, we can obtain
\begin{equation}\label{eqn:b4}
\aligned
&\liminf\limits_{\e\rightarrow0}\e^{-(N+1)}J_\e(\bar{u}_\e)\\
=&\liminf\limits_{\e\rightarrow0}\e^{-(N+1)}\int_{\Omega^1}\frac{1}{|x|^\al}(\frac{1}{2}f(\bar{u}_\e)\bar{u}_\e-F(\bar{u}_\e))dx\\
\geq&\liminf\limits_{\e\rightarrow0}\e^{-(N+1)}\int_{B_{r_\e}(P_\e)\cup B_{r_\e}(P'_\e)\cup B_{r_\e}(Q_\e)}\frac{1}{|x|^\al}(\frac{1}{2}f(\bar{u}_\e)\bar{u}_\e-F(\bar{u}_\e))dx\\
=&\liminf\limits_{\e\rightarrow0}\bigg[\int_{B_{r_\e/\e}(P_\e)}\frac{1}{|\e x+P_{\e}|^\al}(\frac{1}{2}f(\bar{v}_{\e,1})\bar{v}_{\e,1}-F(\bar{v}_{\e,1}))dx\bigg]\\
&+\liminf\limits_{\e\rightarrow0}\bigg[\int_{B_{r_\e/\e}(P'_\e)}\frac{1}{|\e x+P'_{\e}|^\al}(\frac{1}{2}f(\bar{v}_{\e,2})\bar{v}_{\e,2}-F(\bar{v}_{\e,2}))dx\bigg]\\
&+\liminf\limits_{\e\rightarrow0}\bigg[\int_{B_{r_\e/\e}(Q_\e)}\frac{1}{|\e x+Q_{\e}|^\al}(\frac{1}{2}f(\bar{v}_{\e,3})\bar{v}_{\e,3}-F(\bar{v}_{\e,3}))dx\bigg]\\
\geq&\bigg[\int_{\R^{N+1}}\frac{1}{|\bar{P}|^\al}(\frac{1}{2}f(\bar{v}_{1})\bar{v}_{1}-F(\bar{v}_{1}))dx\bigg]
+\bigg[\int_{\R^{N+1}}\frac{1}{|\bar{P}'|^\al}(\frac{1}{2}f(\bar{v}_{2})\bar{v}_{2}-F(\bar{v}_{2}))dx\bigg]\\
&+\bigg[\int_{\R^{N+1}}\frac{1}{|\bar{Q}|^\al}(\frac{1}{2}f(\bar{v}_{3})\bar{v}_{3}-F(\bar{v}_{3}))dx\bigg]\\
\geq&\bigg(|\bar{P}|^{\frac{\al}{2}(N-1)}+|\bar{Q}|^{\frac{\al}{2}(N-1)}+|\bar{P}'|^{\frac{\al}{2}(N-1)}\bigg)
J(U),
\endaligned
\end{equation}
where $U$ is the unique positive solution of equation (\ref{eqn:spacefq}).
Lemma \ref{Lem:limit} implies a contradiction with (\ref{eqn:b4}).
Similarly, we can also show
 that $\bar{u}_\e$ has only one local minimum.
\qed

\begin{remark}\label{r3.8}
Since $u_\e$ has only two
positive local maximums and only two negative local minimums, points $\pm P_\e,\pm Q_\e$ must belong to the $x_{N+1}$-axis for sufficiently small $\e$, that is,
$$\{P_\e,Q_\e\}\subset \mathcal{L}:=\{x\in \Omega:\,\,x=(\mathbf{0},\xi),\,\xi\in\R\}.$$
\end{remark}

 In the following, we will prove the exponential decay of $u_\e$.
 \begin{lemma}\label{Lem:loca-decay}
For any $\delta\in(0,1)$, there exist
$C,c>0$ independent of $\e$ such that for any $x\in \Omega$,
$$
|u_\e(x)|+\e|\nabla u_\e(x)|\leq Ce^{-c(1-\delta)\min\{|x-P_\e|,|x-Q_\e|\}/\e}.
$$
\end{lemma}
 \Proof
 Assume that $\bar{v}_\e^+(y):=\bar{u}_\e^+(\e y+P_\e)$ and $\bar{v}_\e^-(y):=\bar{u}_\e^-(\e y+Q_\e)$.
 Then recalling Lemma \ref{Lem:limit-blow},
 for any $R>0$ there exists constants $\e_0>0$ such that for all
 $\e\in(0,\e_0)$, we have
 $$
 B_{R\e}(P_\e)\subset\Omega_{P_\e}^1\quad\quad B_{R\e}(Q_\e)\subset\Omega_{Q_\e}^1.
 $$
By the elliptic regularity theory and the boundedness of $\{u_\e\}$ in $H_{\sharp}(\Omega)$,
we obtain
$\bar{v}^\pm_\e\rightarrow \bar{v}^\pm$ in $C_{loc}^2(\R^{N+1})$ as $\e\rightarrow 0^+$.
 Moreover,
the standard elliptic $L^q$-estimate with $q>1$ yields $\bar{v}^\pm\in C_{loc}^2(\R^{N+1})\cap W^{2,q}(\R^{N+1})$.
Arguing similarly as in Lemma \ref{Lem:Unique}, we get
$\bar{v}^\pm_\e\rightarrow \bar{v}^\pm$ in $H^1(\R^{N+1})$ as $\e\rightarrow 0^+$. By the standard elliptic estimates,
there exists $C > 0$ (independent $\e$) such that
for any $B_2(y)\in \Omega_{P_\e}^1\cap \Omega_{Q_\e}^1$,
$$
\sup\limits_{B_1(y)}|\bar{v}^\pm_\e|\leq C\|\bar{v}^\pm_\e\|_{L^2(B_2(y))}.
$$
It then follows that
$$
\lim\limits_{|x|\rightarrow+\infty}|\bar{v}^\pm_\e(x)|=0\quad\text{uniformly\,for}\,\e>0\,\text{small\,enough.}
$$
 Thus, for any $\eta>0$, take $R$ large sufficiently, we have
$|\bar{v}^\pm_\e(x)|\leq \eta$ for $|x|\geq R$.
From ($f_1$), we can choose $\eta$ such that for any $\delta>0$, $f(t)/t\leq \delta$ for $|t|\leq \eta$.
Then $v_\e^{\pm}$ satisfies
\begin{equation}\label{eqn:comparison}
\aligned
\left\{
  \begin{array}{ll}
  \displaystyle -\triangle |\bar{v}_\e^{\pm}|+(1-\delta)\frac{|\bar{v}_\e^{\pm}|}{|\e x+P_\e^i|^\al}\leq0&  \,\, |x|\geq R, i=1,2,\\
   | \bar{v}_\e^{\pm}|\leq \eta,& \,\,|x|=R,\\
  \end{array}
\right.
\endaligned
\end{equation}
where $P_\e^1=P_\e$ and $P_\e^2=Q_\e$.
By the maximum principle in a standard way, we show that
$$
|\bar{v}_\e^{\pm}(y)|\leq Ce^{-c(1-\delta)|y|},\quad y\in\R^{N+1}
$$
for $\delta>0$ small enough and some positive constant $C,c$ independently of $\e$.
By a scaling technique, we have
$$
|\bar{u}_\e^{\pm}(x)|\leq Ce^{\frac{-c(1-\delta)|x-P_\e^i|}{\e}},\quad x\in\R^{N+1},\,i=1,2.
$$
By the Harnack inequality, we also have
$$
\bar{u}_\e^{\pm}(x)+\e|\nabla\bar{u}_\e^{\pm}(x)|\leq Ce^{\frac{-c(1-\delta)|x-P_\e^i|}{\e}},\quad x\in\R^{N+1}.
$$
The conclusion of lemma follows from $u_\e\in H_{\sharp}(\Omega)$. The proof is complete.
\qed

\s{Proof of Theorem \ref{Thm:inner-out}}
\vskip0.1in
In this section, we will divide two cases to complete the proof of Theorem \ref{Thm:inner-out}.
\subsection{Case $\eta<2$}
 Note that $\eta<2$, then $\al>0$.
Since $u_\e$ is a least energy sign-changing solution of equation (\ref{eqn:redu-qs}),
it is natural to expect that the points $\{\pm P_\e\},\{\pm Q_\e\}$ should converge to points $\pm\bar{P},\pm\bar{Q}$ in the annulus which have the smallest distance
from the origin. They are indeed the points on the inner boundary.

\Proof
We proceed the proof by contradiction.
In virtue of Lemma \ref{Lem:axis}, without loss of generality, we assume that $\{P_\e\},\{Q_\e\}\subset \Omega^1$,
and that $P_\e=(\mathbf{0},\frac{a^2}{2}+\nu_\e)$ and $P_\e$ converge to point $\bar{P}$ with
$|\bar{P}|=\frac{a^2}{2}+\nu$ for some $\nu>0$. That is to say, $\nu_\e\rightarrow\nu$ as $\e\rightarrow0$.
Without loss of generality, we also assume
$$
Q_\e=\big(\mathbf{0},\frac{a^2}{2}+\bar{\nu}_\e\big)\quad\quad\text{and}\quad\quad \bar{\nu}_\e<\nu_\e.
$$
Now we divide two cases to state our proof.\\
\emph{Case 1:} $\bar{\nu}_\e\rightarrow0$, that is, $|\bar{Q}|=\frac{a^2}{2}$ and $|\bar{P}-\bar{Q}|=\nu$.
Consider the ball $B_{\frac{\nu}{4}}(\hat{P})$
with center at $\hat{P}=(\mathbf{0},\frac{a^2+\nu}{2})$.
Set a cut off function $\varphi\in C_0^\infty(\R^{N+1})$ with $\varphi(x)\equiv1$ for $x\in B_{\frac{\nu}{8}}(0)$ and
$\varphi(x)\equiv0$ for $x\in \R^{N+1}\setminus{B_{\frac{\nu}{4}}(0)}$.
Define the function
$$
h^+_\e(x)=\varphi(x-\hat{P})w(\frac{x-\hat{P}}{\e})
+\varphi(x+\hat{P})w(\frac{x+\hat{P}}{\e})=:\mathbf{A_\e}+\mathbf{B_\e},
$$
where $w$ is the unique positive solution of equation (\ref{eqn:redu-qs}) with $d=|\hat{P}|$.
It is easy to see that there exists a unique $t_\e>0$ such that $J'_\e(t_\e h^+_\e)t_\e h^+_\e=0$.
We now show $t_\e\rightarrow1$ as $\e\rightarrow0^+$. Indeed, based on the definition of $h^+_\e(x)$, the supports of
$\mathbf{A_\e}$ and $\mathbf{B_\e}$ are disjoint each other. So we have
\begin{equation}\label{eqn:guji6}
\aligned
&J'_\e(t_\e \mathbf{A_\e})t_\e \mathbf{A_\e}\\
=&\e^2t_\e^2\int_{\Omega}|\nabla  \mathbf{A_\e}|^2dx+t_\e^2
\int_{\Omega}\frac{\mathbf{A_\e}^2}{|x|^\alpha}dx-\int_{\Omega}\frac{f(t_\e \mathbf{A_\e})t_\e \mathbf{A_\e}}{|x|^\alpha}dx\\
=&\e^{N+1}t_\e^2\bigg[\int_{B_{\frac{\nu}{4\e}}(\hat{P})}|\nabla (\varphi(\e x)w)|^2dx+
\int_{B_{\frac{\nu}{4\e}}(\hat{P})}\frac{|\varphi(\e x)w|^2}{|\e x+\hat{P}|^\alpha}dx\\
&-\int_{B_{\frac{\nu}{4\e}}(\hat{P})}\frac{f(t_\e \varphi(\e x)w) \varphi(\e x)w}{t_\e|\e x+\hat{P}|^\alpha}dx\bigg]
=0,
\endaligned
\end{equation}
which implies by  the fact that $|w(x)|\leq C\exp(-\frac{\delta}{|\hat{P}|^\al}|x|)$ for some $\delta>0$ that
\begin{equation}\label{eqn:guji7}
\aligned
&\int_{\R^{N+1}}|\nabla w|^2dx+
\int_{\R^{N+1}}\frac{|w|^2}{|\hat{P}|^\alpha}dx+O(\e)\\
&=\int_{B_{\frac{\nu}{4\e}}(\hat{P})}|\nabla (\varphi(\e x)w)|^2dx+
\int_{B_{\frac{\nu}{4\e}}(\hat{P})}\frac{|\varphi(\e x)w|^2}{|\e x+\hat{P}|^\alpha}dx\\
&\geq\int_{B_{\frac{\nu}{8\e}}(\hat{P})}\frac{f(t_\e w) w^2}{t_\e w|\e x+\hat{P}|^\alpha}dx.
\endaligned
\end{equation}
By  ($f_2$), we deduce that $\{t_\e\}$ is bounded in $\R$. Without loss of generality, we assume
$t_\e\rightarrow t_0\geq0$ as $\e\rightarrow0$.
Moreover, from (\ref{eqn:guji6}) we deduce that
\begin{equation}\label{eqn:guji8}
\int_{B_{\frac{\nu}{8\e}}(\hat{P})}|\nabla w|^2dx+
\int_{B_{\frac{\nu}{8\e}}(\hat{P})}\frac{|w|^2}{|\e x+\hat{P}|^\alpha}dx
<\int_{B_{\frac{\nu}{4\e}}(\hat{P})}\frac{f(t_\e w)w}{t_\e|\e x+\hat{P}|^\alpha}dx.
\end{equation}
It then follows that $t_0\not=0$. In virtue of (\ref{eqn:guji6}), passing the limit as $\e\rightarrow0^+$, we get
\begin{equation}\label{eqn:guji9}
\int_{\R^{N+1}}|\nabla w|^2dx+
\int_{\R^{N+1}}\frac{|w|^2}{|\hat{P}|^\alpha}dx
-\int_{\R^{N+1}}\frac{f(t_0w)w^2}{t_0w|\hat{P}|^\alpha}dx=0.
\end{equation}
Recalling that $w$ is the unique positive solution of equation (\ref{eqn:redu-qs}) with $d=|\hat{P}|$,
Remark 1.1 yields $t_0=1$. Thus, using (\ref{eqn:Uex}) and the dominate convergence theorem yields
\begin{equation}\label{eqn:guji10}
\aligned
&\e^{-(N+1)}J_\e(t_\e h^+_\e)\\
=&\frac{t_\e^2\e^{(1-N)}}{2}\int_{\Omega}|\nabla h^+_\e|^2dx+\frac{t_\e^2\e^{-(N+1)}}{2}
\int_{\Omega}\frac{|h^+_\e|^2}{|x|^\alpha}dx-\e^{-(N+1)}\int_{\Omega}\frac{F(t_\e h^+_\e)}{|x|^\alpha}dx\\
=&t_\e^2\bigg[\frac{1}{2}\int_{B_{\frac{\nu}{4\e}}(\hat{P})}|\nabla (\varphi(\e x)w)|^2dx+\frac{1}{2}
\int_{B_{\frac{\nu}{4\e}}(\hat{P})}\frac{|\varphi(\e x)w|^2}{|\e x+\hat{P}|^\alpha}dx
-\int_{B_{\frac{\nu}{4\e}}(\hat{P})}\frac{F(t_\e \varphi(\e x)w)}{t_\e^2|\e x+\hat{P}|^\alpha}dx\bigg]\\
&+t_\e^2\bigg[\frac{1}{2}\int_{B_{\frac{\nu}{4\e}}(-\hat{P})}|\nabla (\varphi(\e x)w)|^2dx+\frac{1}{2}
\int_{B_{\frac{\nu}{4\e}}(-\hat{P})}\frac{|\varphi(\e x)w|^2}{|\e x-\hat{P}|^\alpha}dx
-\int_{B_{\frac{\nu}{4\e}}(-\hat{P})}\frac{F(t_\e \varphi(\e x)w)}{t_\e^2|\e x-\hat{P}|^\alpha}dx\bigg]\\
=& 2\bigg[\frac{1}{2}\int_{\R^{N+1}}|\nabla w|^2dx+\frac{1}{2}
\int_{\R^{N+1}}\frac{w^2}{|\hat{P}|^\alpha}dx
-\int_{\R^{N+1}}\frac{F(w)}{|\hat{P}|^\alpha}dx+O(\e)\bigg].
\endaligned
\end{equation}
Thus, by equation (\ref{eqn:redu-qs}) and $\alpha=1-\frac\eta2>0$, we have that
\begin{equation}\label{eqn:guji11}
\lim\limits_{\e\rightarrow0}\e^{-(N+1)}J_\e(t_\e h^+_\e)=2|\hat{P}|^{\frac{(N-1)\al}{2}}J(U)
<2|\bar{P}|^{\frac{(N-1)\al}{2}}J(U).
\end{equation}
Define
$$
h^-_\e(x):=-\psi(x-Q_\e)u_\e(\frac{x-Q_\e}{\e })-\psi(x+Q_\e)u_\e(\frac{x+Q_\e}{\e }),
$$
 where $\psi$ is a non-negative
smooth radial function supported in $B_{2r_\e}(0)$ with $|\nabla \phi|\leq \frac{2}{r_\e}$ and
$$
\phi(r)=\left\{
  \begin{array}{ll}
  \displaystyle  1,&  \mbox{for}\,\, r\in[0,r_\e],\\
   0,&  \mbox{for}\,\, r\in[2r_\e,+\infty),
  \end{array}
\right.
$$
where $r_\e$ is chosen so that $4r_\e=\min\{dist(Q_\e,\partial\Omega^1),dist(Q_\e,\Theta)\}$. Note that
there exists a unique $s_\e>0$ such that $J'_\e(s_\e h^-_\e)s_\e h^-_\e=0$.
In virtue of Lemma \ref{Lem:loca-decay}, using the similar argument as above, we can also obtain
$s_\e\rightarrow1$ and furthermore
\begin{equation}\label{eqn:guji12}
\lim\limits_{\e\rightarrow0}\e^{-(N+1)}J_\e(s_\e h^-_\e)=2|\bar{Q}|^{\frac{(N-1)\al}{2}}J(U).
\end{equation}
Note that $h^+_\e$ and $h^-_\e$ have disjoint supports for $\e$ small. Let us define
$$
h_\e=t_\e h^+_\e+s_\e h^-_\e.
$$
Then $h_\e\in\mathcal{N}_\e$ and by (\ref{eqn:guji11}) and (\ref{eqn:guji12}), we have
$$
\aligned
\lim\limits_{\e\rightarrow0}\e^{-(N+1)}c_\e&\leq \lim\limits_{\e\rightarrow0}\e^{-(N+1)}J_\e(h_\e)\\
&=\lim\limits_{\e\rightarrow0}\bigg[\e^{-(N+1)}J_\e(t_\e h^+_\e)+\e^{-(N+1)}J_\e(s_\e h^-_\e)\bigg]\\
&<2(|\bar{P}|^{\frac{(N-1)\al}{2}}+|\bar{Q}|^{\frac{(N-1)\al}{2}})J(U),
\endaligned
$$
which contradicts with Lemma \ref{Lem:limit}. \\
\emph{Case 2:}  $\bar{\nu}_\e\rightarrow \bar{\nu}>0$, that is, $|\bar{Q}|=\frac{a^2}{2}+\bar{\nu}$
and $|\bar{P}-\bar{Q}|=\nu-\bar{\nu}\geq0$. Consider the ball $B_{\frac{\bar{\nu}}{4}}(\hat{Q})$
with center at $\hat{Q}=(\mathbf{0},\frac{a^2+\bar{\nu}}{2})$.
Set one radial cutoff function $\varphi\in C_0^\infty(\R^{N+1})$ with $\varphi(x)\equiv1$ for $x\in B_{\frac{\bar{\nu}}{8}}(0)$ and
$\varphi(x)\equiv0$ for $x\in \R^{N+1}\setminus{B_{\frac{\bar{\nu}}{4}}(0)}$.
Define the function
$$
g^-_\e(x)=-\varphi(x-\hat{Q})w(\frac{x-\hat{Q}}{\e})-\varphi(x+\hat{Q})w(\frac{x+\hat{Q}}{\e}),
$$
where $w$ is the unique positive solution of equation (\ref{eqn:redu-qs}) with $d=|\hat{Q}|$.
There also exists $\bar{t}_\e>0$ such that $J'_\e(\bar{t}_\e g^-_\e)\bar{t}_\e g^-_\e=0$.
Arguing as in \emph{Case 1},
we can obtain the similar estimate as (\ref{eqn:guji11})
\begin{equation}\label{eqn:guji13}
\lim\limits_{\e\rightarrow0}\e^{-(N+1)}J_\e(\bar{t}_\e g^-_\e)=2|\hat{Q}|^{\frac{(N-1)\al}{2}}J(U)
<2|\bar{Q}|^{\frac{(N-1)\al}{2}}J(U),
\end{equation}
where $\bar{t}_\e \rightarrow1$ as $\e\rightarrow0^+$ and $\al>0$. Define
$$
g^+_\e(x):=\psi(x-P_\e)u_\e(\frac{x-P_\e}{\e })+\psi(x+P_\e)u_\e(\frac{x+P_\e}{\e }),
$$
where $\psi$ is a radial smooth cut off function supported in $B_{2r_\e}(0)$ with $|\nabla \phi|\leq \frac{2}{r_\e}$ and
$$
\phi(r)=\left\{
  \begin{array}{ll}
  \displaystyle 1,&  \mbox{for}\,\, r\in[0,r_\e],\\
   0,&  \mbox{for}\,\, r\in[2r_\e,+\infty),
  \end{array}
\right.
$$
where $4r_\e=\min\{dist(P_\e,\partial\Omega^1),dist(P_\e,\Theta)\}$. It follows from Lemma \ref{Lem:asym} that $r_\e/\e\rightarrow+\infty$. Using the similar argument as \emph{Case 1},
we have $J'_\e(\bar{s}_\e g^+_\e)\bar{s}_\e g^+_\e=0$ for some $\bar{s}_\e>0$, and
\begin{equation}\label{eqn:guji14}
\lim\limits_{\e\rightarrow0}\e^{-(N+1)}J_\e(s_\e g^+_\e)=2|\bar{P}|^{\frac{(N-1)\al}{2}}J(U),
\end{equation}
where $\bar{s}_\e \rightarrow1$ as $\e\rightarrow0^+$. Hence,
$$
g_\e(x)=\bar{s}_\e g^+_\e+\bar{t}_\e g^-_\e\in\mathcal{N}_\e.
$$
Combining (\ref{eqn:guji13}) and (\ref{eqn:guji14}), we have
$$
\aligned
\lim\limits_{\e\rightarrow0}\e^{-(N+1)}c_\e\leq& \lim\limits_{\e\rightarrow0}\e^{-(N+1)}J_\e(g_\e)\\
=&\lim\limits_{\e\rightarrow0}\bigg[\e^{-(N+1)}J_\e(\bar{s}_\e g^+_\e)+\e^{-(N+1)}J_\e(\bar{t}_\e g^-_\e)\bigg]\\
<&2(|\bar{Q}|^{\frac{(N-1)\al}{2}}+|\bar{P}|^{\frac{(N-1)\al}{2}})J(U),
\endaligned
$$
which contradicts with Lemma \ref{Lem:limit}. As a consequence, $|\bar{P}|=|\bar{Q}|=\frac{a^2}{2}$.
Recalling Lemma \ref{Lem:mox}, Theorem \ref{Thm:exist} and Lemma \ref{Lem:loca-decay}, equation (\ref{eqn:redu-qs-})
has a nonradial nodal solution concentrating at four points $\pm\bar{P},\pm\bar{Q}$ in $\Omega$. More precisely,
$$
\pm\bar{P}=(\mathbf{0},\pm \frac{a^2}{2})\quad \pm\bar{Q}=(\mathbf{0},\pm \frac{a^2}{2}).
$$
That is, $\bar{P}=\bar{Q}$.
Then the corresponding solution of equation
(\ref{eqn:qs}), still denoted by $u_\e$, concentrates exactly at two orthogonal $(N-1)$-dimensional spheres in surface $|x|=a$, placed
the angle $\theta=0$ and $\theta=\frac{\pi}{2}$.
\qed
\subsection{Case $\eta>2$} Contrary to Case $\eta<2$, we expect that the points $\{\pm P_\e\},\{\pm Q_\e\}$
should converge to points $\pm\bar{P},\pm\bar{Q}$ in the annulus which have the largest distance
from the origin, since $\al=1-\frac\eta2<0$. There are indeed the points on the outer boundary.

\Proof
We begin to prove by contradiction.
According to Lemma \ref{Lem:axis},
without loss of generality, we assume that $\{P_\e\},\{Q_\e\}\subset \Omega^1$,
and that $P_\e=(\mathbf{0},\frac{b^2}{2}-\nu_\e)$ and $P_\e$ converge to point $\bar{P}$ with
$|\bar{P}|=\frac{b^2}{2}-\nu$ for some $\nu>0$. That is to say, $\nu_\e\rightarrow\nu$.
Without loss of generality, we also assume
$$
Q_\e=\big(\mathbf{0},\frac{b^2}{2}-\bar{\nu}_\e\big)\quad\quad\text{and}\quad\quad \bar{\nu}_\e<\nu_\e.
$$
If $\bar{\nu}_\e\rightarrow0$, that is, $|\bar{Q}|=\frac{b^2}{2}$, similarly to Case $\eta<2$, we have

$$
\aligned
\lim\limits_{\e\rightarrow0}\e^{-(N+1)}c_\e&\leq 2(|\hat{P}|^{\frac{(N-1)\al}{2}}+|\bar{Q}|^{\frac{(N-1)\al}{2}})J(U)\\
&<2(|\bar{P}|^{\frac{(N-1)\al}{2}}+|\bar{Q}|^{\frac{(N-1)\al}{2}})J(U),
\endaligned
$$
where $\hat{P}=(\mathbf{0},\frac{b^2-\nu}{2})$ and $\alpha=1-\frac\eta2<0$ is used. This is a contradiction.  If $\bar{\nu}_\e\rightarrow \bar{\nu}>0$, that is, $|\bar{Q}|=\frac{b^2}{2}-\bar{\nu}$, similarly to Case $\eta<2$, we have

$$
\aligned
\lim\limits_{\e\rightarrow0}\e^{-(N+1)}c_\e&\leq 2(|\bar{P}|^{\frac{(N-1)\al}{2}}+|\hat{Q}|^{\frac{(N-1)\al}{2}})J(U)\\
&<2(|\bar{P}|^{\frac{(N-1)\al}{2}}+|\bar{Q}|^{\frac{(N-1)\al}{2}})J(U),
\endaligned
$$
where $\hat{Q}=(\mathbf{0},\frac{b^2-\bar{\nu}}{2})$ and $\alpha=1-\frac\eta2<0$ is used. This is a contradiction again. Thus, we get that $|\bar{P}|=|\bar{Q}|=\frac{b^2}{2}$. Therefore,
equation (\ref{eqn:redu-qs-})
has a nonradial nodal solution concentrating at four points $\pm\bar{P},\pm\bar{Q}$
in $\Omega$.  More precisely,
$$
\pm\bar{P}=(\mathbf{0},\pm \frac{b^2}{2})\quad \pm\bar{Q}=(\mathbf{0},\pm \frac{b^2}{2}).
$$
That is, $\bar{P}=\bar{Q}$.
The corresponding solution of equation
(\ref{eqn:qs}), still denoted by $u_\e$, concentrates exactly at two orthogonal $(N-1)$-dimensional spheres in surface $|x|=b$, placed
the angle $\theta=0$ and $\theta=\frac{\pi}{2}$.
\qed

\s{Proof of Theorem \ref{Thm:middle}}
\vskip0.1in
In this section, we consider the case $\eta=2$, that is, $\al=0$.
Our aim is to investigate the location of the spikes of concentrating solution as $\e\rightarrow0^+$.
Define the function
$$
\mathcal{E}_{a_1,a_2}(P_1,P_2):=\min\bigg\{\frac{a_1}{a_1+a_2}|P_1-P_2|,\frac{a_2}{a_1+a_2}|P_1-P_2|,
dist(P_1,\partial \Omega)\},dist(P_2,\partial \Omega)\}\bigg\},
$$
where $P_1,P_2\in\mathcal{L}^1:=\mathcal{L}\cap \Omega^1$ and $a_1,a_2>0$. Here, $\mathcal{L}$ has been defined in Remark \ref{r3.8}. Let
$$
\mathcal{F}_{a_1,a_2}=\max\limits_{(P_1,P_2)\in \mathcal{L}^1\times\mathcal{L}^1}\mathcal{E}_{a_1,a_2}(P_1,P_2).
$$
We now state the following upper-bound of $c_\e$ defined in Theorem \ref{Thm:exist}.
\begin{lemma}\label{Lem:2upper}
For $\delta>0$ small enough, we have
$$
c_\e\leq 4\e^{N+1}(J(U)+Ce^{-\frac{2(1-\delta)\mathcal{F}_{a_1,a_2}}{\e}})
$$
\end{lemma}
\Proof
Let us define
$$
\Theta_{a_1,a_2}(P_1,P_2)=\{x\in\Omega|\,\,a_1|x-P_1|=a_2|x-P_2|\}
$$
for every $(P_1,P_2)\in \mathcal{L}^1\times\mathcal{L}^1$, and set
$$
R_i=\min\bigg\{dist(P_i,\partial\Omega),dist(P_i,\Theta_{a_1,a_2}(P_1,P_2))\bigg\}
$$
respectively for $i=1,2$.
It is easy to see from the definition of $\Theta_{a_1,a_2}(P_1,P_2)$ that $B_{R_1}(P_1)\cap B_{R_2}(P_2)=\emptyset$
and $B_{R_i}(P_i)\subset\Omega$ for $i=1,2$. Moreover,
$$
\aligned
&dist(P_1,\Theta_{a_1,a_2}(P_1,P_2))=\frac{a_2}{a_1+a_2}|P_1-P_2|,\\
&dist(P_2,\Theta_{a_1,a_2}(P_1,P_2))=\frac{a_1}{a_1+a_2}|P_1-P_2|,
\endaligned
$$
which implies that $R_i\geq \mathcal{E}_{a_1,a_2}(P_1,P_2)$ for $i=1,2$.
 Let $\psi^i_\e$ be smooth radial symmetric and decreasing functions so that
$0\leq\psi_\e^i\leq1$ and $\psi_\e^i(x)\equiv1$ for $|x|\leq \frac{R_i}{\e}-1$ and $\psi_\e^i(x)\equiv0$
for $|x|\geq \frac{R_i}{\e}$. Then it follows that
$$
U_\e^i(x):=U(\frac{x-P_i}{\e})\psi^i_\e(\frac{x-P_i}{\e})+U(\frac{x+P_i}{\e})\psi^i_\e(\frac{x+P_i}{\e}),\,\, i=1,2,
$$
where $U$ is the unique positive solution of equation (\ref{eqn:spacefq}).
Obviously,  $U_\e^i\in H_{\sharp}(\Omega)$.
By ($f_1$)-($f_3$) and Theorem \ref{Thm:exist}, it is easy to see that there exist $t_\e^i>0,\,i=1,2$ such that
$t_\e^1U_\e^1-t_\e^2U_\e^2\in \mathcal{N}_\e$.
Using the similar arguments as in Theorem \ref{Thm:inner-out}, we can obtain that $\{t_\e^i\}$
is bounded uniformly for $\e$ small.
 Recalling (\ref{eqn:Uex}), we have
\begin{equation}\label{eqn:2-1}
U(x)+|\nabla U(x)|\leq Ce^{-(1-\delta)|x|}
\end{equation}
for $\delta>0$ small. Since the supports of $U(\frac{x-P_i}{\e})\psi^i_\e(\frac{x-P_i}{\e})$ and
$U(\frac{x+P_i}{\e})\psi^i_\e(\frac{x+P_i}{\e})$
are disjoint, it follows from (\ref{eqn:2-1}) and the boundedness of $\{t_\e^i\}$ for small $\e$ that
\begin{equation}\label{eqn:2-2}
\aligned
&J_\e(t_\e^1U_\e^1)\\
=&J_\e(t_\e^1U(\frac{x-P_1}{\e})\psi^1_\e(\frac{x-P_1}{\e}))+J_\e(t_\e^1U(\frac{x-P_1}{\e})\psi^1_\e(\frac{x-P_i}{\e}))\\
=&\frac{(t_\e^1)^2\e^2}{2}\int_{B_{R_1}(P_1)}|\nabla (U(\frac{x-P_1}{\e})\psi^1_\e(\frac{x-P_1}{\e}))|^2dx+\frac{(t_\e^1)^2}{2}
\int_{B_{R_1}(P_1)}|U(\frac{x}{\e}-P_1)\psi^1_\e(\frac{x-P_1}{\e})|^2dx\\
&-\int_{B_{R_1}(P_1)}F(t_\e^1 U(\frac{x-P_1}{\e})\psi^1_\e(\frac{x-P_1}{\e}))dx
+\frac{(t_\e^1)^2\e^2}{2}\int_{B_{R_1}(-P_1)}|\nabla (U(\frac{x+P_1}{\e})\psi^1_\e(\frac{x+P_1}{\e}))|^2dx\\
&+\frac{(t_\e^1)^2}{2}
\int_{B_{R_1}(-P_1)}|U(\frac{x+P_1}{\e})\psi^1_\e(\frac{x+P_1}{\e})|^2dx
-\int_{B_{R_1}(-P_1)}F(t_\e^1 U(\frac{x+P_1}{\e})\psi^1_\e(\frac{x+P_1}{\e}))dx\\
=&2\e^{N+1}\bigg[\frac{(t_\e^1)^2}{2}\int_{B_{\frac{R_1}{\e}}(0)}\bigg(|\nabla (U(x)\psi_\e^1(x))|^2+|U(x)\psi_\e^1( x)|^2\bigg)dx
-\int_{B_{\frac{R_1}{\e}}(0)}F(t_\e^1 U(x)\psi_\e^1(x))dx\bigg]\\
\leq& 2\e^{N+1}\bigg[\max\limits_{t\in(0,\infty)}J(tU)+Ce^{-2(1-\delta)(\frac{R_1}{\e}-1)}\bigg]\\
\leq& 2\e^{N+1}\bigg[J(U)+Ce^{-2(1-\delta)\frac{R_1}{\e}}\bigg].
\endaligned
\end{equation}
Using the almost same argument as in (\ref{eqn:2-2}), we have
\begin{equation}\label{eqn:2-3}
J_\e(t_\e^2U_\e^2)
\leq 2\e^{N+1}\bigg[J(U)+Ce^{-2(1-\delta)\frac{R_2}{\e}}\bigg].
\end{equation}
Combining (\ref{eqn:2-2}) and (\ref{eqn:2-3}), using the fact that
the supports of $U_\e^1$ and $U_\e^2$ are disjoint,
we have
\begin{equation}\label{eqn:2-4}
\aligned
c_\e&=J_\e(u_\e)\leq J_\e(t_\e^1U_\e^1-t_\e^2U_\e^2)\\
&=J_\e(t_\e^1U_\e^1)+J_\e(t_\e^2U_\e^2)\\
&\leq 4\e^{N+1}\bigg[J(U)+Ce^{-2(1-\delta)\frac{\min\{R_1,R_2\}}{\e}}\bigg]\\
&\leq 4\e^{N+1}\bigg[J(U)+Ce^{-2(1-\delta)\frac{\mathcal{E}_{a_1,a_2}(P_1,P_2)}{\e}}\bigg].
\endaligned
\end{equation}
The conclusion follows from the arbitrariness of $P_1,P_2$.
The proof is complete.
\qed

Recall that $u_\e \in H_{\sharp}(\Omega)$ is a nodal solution of equation (\ref{eqn:redu-qs}). Let us define
$$
R_{\e,+}:=\max\{R>0|\,B_{R}(P_\e)\subset \Omega_{\e,+}^1\},\quad R_{\e,-}:=\max\{R>0|\,B_{R}(Q_\e)\subset \Omega_{\e,-}^1\},
$$
where $\Omega_{\e,\pm}^1$ are the support sets of $\bar{u}_\e^\pm$. In virtue of Lemma \ref{Lem:limit-blow}, we have
$\frac{R_{\e,\pm}}{\e}\rightarrow\infty$ as $\e\rightarrow0^+$.
 Obviously,
$$
B_{R_{\e,+}}(P_\e)\cap B_{R_{\e,-}}(Q_\e)=\emptyset,\quad\text{ and}\quad
B_{R_{\e,+}}(P_\e),B_{R_{\e,-}}(Q_\e)\subset \Omega^1\subset\Omega.
$$
Set
$$
\aligned
&v_\e^+(y)=v_{1,\e}^+(y)+v_{2,\e}^+(y):=\bar{u}_\e^+(\e y+P_\e)\psi(\frac{|y|}{R_{\e,+}})+\tilde{u}_\e^+(\e y-P_\e)\psi(\frac{|y|}{R_{\e,+}}),\\
&v_\e^-(y)=v_{1,\e}^-(y)+v_{2,\e}^-(y):=\bar{u}_\e^-(\e y+Q_\e)\psi(\frac{|y|}{R_{\e,-}})+\tilde{u}_\e^-(\e y-Q_\e)\psi(\frac{|y|}{R_{\e,-}}),
\endaligned
$$
where $\psi$ is the following smooth radial and decreasing cut-off function
\begin{equation}\label{eqn:2-5-0}
\psi(r)=\left\{
  \begin{array}{ll}
  \displaystyle  1&  \mbox{for}\,\, r\in[0,\frac{1}{\e}-\delta'],\\
   0,&  \mbox{for}\,\, r\in[\frac{1}{\e},+\infty),
  \end{array}
\right.
\end{equation}
with $|\psi'(r)|\leq C$ and $\delta'>0$ small constant. It is easy to see that
$v_\e^{\pm}(y)\in H_0^1(B_{\frac{R_{\e,\pm}}{\e}})$. Here, $\bar{u}_\e,\tilde{u}_\e$ have been defined in Lemma \ref{Lem:asym}.
Recalling Theorem \ref{Thm:exist} and Lemma \ref{Lem:loca-decay}, we have for any $t,s\in\R$
\begin{equation}\label{eqn:2-5}
\aligned
J_\e(u_\e)\geq& J_\e(tu_\e^+)+J_\e(su_\e^-)\\
=&J_{\e,\Omega^1}(t\bar{u}_\e^+)+J_{\e,\Omega^2}(t\tilde{u}_\e^+)+J_{\e,\Omega^1}(s\bar{u}_\e^-)+J_{\e,\Omega^2}(s\tilde{u}_\e^-)\\
=& \epsilon^{(N+1)}\bigg[J_{\Omega^1_{P_\e,+}}(t\bar{u}_\e^+(\e y+P_\e))+J_{\Omega^2_{P_\e,-}}(t\tilde{u}_\e^+(\e y-P_\e))\\
&+J_{\Omega^1_{Q_\e,+}}(s\bar{u}_\e^-(\e y+Q_\e))+J_{\Omega^2_{Q_\e,-}}(s\tilde{u}_\e^-(\e y-Q_\e))\bigg]\\
\geq &\e^{N+1}\bigg[J_{B_{\frac{R_{\e,+}}{\e}}(0)}(tv_{1,\e}^+)+J_{B_{\frac{R_{\e,+}}{\e}}(0)}(tv_{2,\e}^+)
+J_{B_{\frac{R_{\e,-}}{\e}}(0)}(sv_{1,\e}^-)+J_{B_{\frac{R_{\e,-}}{\e}}(0)}(sv_{2,\e}^-)\bigg]\\
&-C(t)e^{-2(1-\delta)R_{\e,+}(\frac{1}{\e}-\delta')}-C(s)e^{-2(1-\delta)R_{\e,-}(\frac{1}{\e}-\delta')}.
\endaligned
\end{equation}
Here, $J_{\e,\Omega^i},$ denote
functional $J_{\e}$ whose integral region are $\Omega^i, i=1,2$. $J_{\Omega^1_{P_\e,+}}$, $J_{\Omega^2_{P_\e,-}}$,
$J_{\Omega^1_{Q_\e,+}}$, $J_{\Omega^2_{Q_\e,-}}$ and $J_{B_{\frac{R_{\e,\pm}}{\e}}(0)}$ denote
functional $J$ whose integral region are $\Omega^1_{P_\e,+}$, $\Omega^2_{P_\e,-}$,
$\Omega^1_{Q_\e,+}$, $\Omega^2_{Q_\e,-}$ and
$B_{\frac{R_{\e,\pm}}{\e}}(0)$, respectively. Moreover, $C(t)$ and $C(s)$ are bounded from above since we can always choose
$t,s$ are bounded above.

It is easy to prove that for $\e>0$ small enough,
\begin{equation}\label{eqn:half-eq1,2}
\left\{
  \begin{array}{ll}
  \displaystyle -\triangle u+u=f(u),&  \,\,\, x\in B_{\frac{R_{\e,\pm}}{\e}}(0),\\
    u=0,&\,\,\,x\in \partial B_{\frac{R_{\e,\pm}}{\e}}(0)
  \end{array}
\right.
\end{equation}
have positive ground state solutions $v_{\e,+}$ and $v_{\e,-}$ in
$H_0^1(B_{\frac{R_{\e,+}}{\e}}(0))$ and $H_0^1(B_{\frac{R_{\e,-}}{\e}}(0))$, respectively. Using the well-known Gidas-Ni-Nirenberg's theorem
\cite{Gidas81}, we see that $v_{\e,+}$ and $v_{\e,-}$ are both radial symmetrical, and that $v_{\e,+},v_{\e,-}$ are nonincreasing with respect to $r$.
It is easy to see from ($f_2$) that $\{v_{\e,+}\}$ and $\{v_{\e,-}\}$ are bounded uniformly for $\e$. Hence,
since $v_{\e,+},v_{\e,-}$ are nonincreasing, we have
\begin{equation}\label{eqn:decay0}
\lim\limits_{|x|\rightarrow\infty}v_{\e,\pm}(x)=0\quad\text{uniformly\,\,for\,\,}\e>0\,\,\text{ small\,\,enough.}
\end{equation}
Then, by the comparison principle, we see that for $\e>0$ small enough,
there exists a constant $c,C>0$ satisfying
\begin{equation}\label{eqn:decayv}
ce^{-(1+\delta)R_{\e,\pm}(\frac{1}{\e}-\delta')}\leq v_{\e,\pm}(R_{\e,\pm}(\frac{1}{\e}-\delta'))\leq Ce^{-(1-\delta)R_{\e,\pm}(\frac{1}{\e}-\delta')}
\end{equation}
with any $\delta\in(0, 1)$ independent of $\e$. As in the proof of Theorem 4.1 in \cite{Byeon05}, we can extend $v_{\e,\pm}$ to the whole space $\R^{N+1}$
and denote it by
$$
w_{\e,\pm}\equiv\left\{
  \begin{array}{ll}
  \displaystyle v_{\e,\pm},&  \,\,\, \text{for}\quad|x|\leq R_{\e,\pm}(\frac{1}{\e}-\delta'), \\
    V_{\e,\pm},&\,\,\,\text{for}\quad |x|\geq R_{\e,\pm}(\frac{1}{\e}-\delta'),
  \end{array}
\right.
$$
where $V_{\e,\pm}$ is a radial symmetric positive solution of the following equation
$$
\left\{
  \begin{array}{ll}
  \displaystyle -\triangle u+u=f(u)&  \,\,\, \text{for}\quad|x|> R_{\e,\pm}(\frac{1}{\e}-\delta'), \\
    u(R_{\e,\pm}(\frac{1}{\e}-\delta'))=v_{\e,\pm}(R_{\e,\pm}(\frac{1}{\e}-\delta')),&\,\,\,\text{and}\quad \lim\limits_{|x|\rightarrow\infty}u(x)=0.
  \end{array}
\right.
$$
Using the comparison principle and the above facts, we have that for any $\delta\in(0, 1)$,
there exists a constant $C>0$ such that
\begin{equation}\label{eqn:decayV}
V_{\e,\pm}(r)+|V'_{\e,\pm}(r)|\leq C e^{-(1-\delta)r}\quad \text{for}\,\, |x|\geq R_{\e,\pm}(\frac{1}{\e}-\delta').
\end{equation}
Thus, there exist exactly $t_{\e,\pm}>0$ such that
\begin{equation}\label{eqn:zhengfu}
J'(t_{\e,\pm}w_{\e,\pm})t_{\e,\pm}w_{\e,\pm}=0.
\end{equation}
Obviously, $J(t_{\e,\pm}w_{\e,\pm})\geq J(U)$, where $U$ is the unique positive solution of equation (\ref{eqn:spacefq}).
Moreover, we state the following estimates.
\begin{lemma}\label{Lem:estimate}
\begin{equation}\label{eqn:zhengfu2}
J_{B_{\frac{R_{\e,\pm}}{\e}}(0)}(t_{\e,\pm}v_{\e,\pm})
\geq J(t_{\e,\pm}w_{\e,\pm}) +Ce^{-2(1+\delta)R_{\e,\pm}(\frac{1}{\e}-\delta')},
\end{equation}
where $C,\delta>0$ is independent of $\e$ and $\delta'$ has been given in (\ref{eqn:2-5-0}).
\end{lemma}
\Proof
 Note that by the standard elliptic estimates,  $w_{\e,\pm}\rightarrow U$ in
$H^1(\R^{N+1})\cap C_{loc}^2(\R^{N+1})$ as $\e\rightarrow0^+$, where $U$ is
the unique positive solution of equation (\ref{eqn:spacefq}).
Thus, it is easy to see from (\ref{eqn:zhengfu}) that
$t_{\e,\pm}\rightarrow1$ as $\e\rightarrow0^+$. And then,
$$
J(t_{\e,\pm}w_{\e,\pm})\rightarrow J(U),\quad\quad J(w_{\e,\pm})\rightarrow J(U),\quad\text{as}\,\,\e\rightarrow0^+.
$$
Let us define $g_{\e,\pm}(t):=J(tw_{\e,\pm})$, then by ($f_2$) we have as $\e\rightarrow0^+$
$$
\aligned
g''_{\e,\pm}(t_{\e,\pm})&=\int_{\R^{N+1}}(|\nabla w_{\e,\pm}|^2+|w_{\e,\pm}|^2)dx-\frac{1}{t^2_{\e,\pm}}\int_{\R^{N+1}}f'(t_{\e,\pm}w_{\e,\pm})|t_{\e,\pm} w_{\e,\pm}|^2dx\\
&\leq\int_{\R^{N+1}}(|\nabla w_{\e,\pm}|^2+|w_{\e,\pm}|^2)dx-\frac{1}{t^2_{\e,\pm}}\int_{\R^{N+1}}\mu f(t_{\e,\pm}w_{\e,\pm})t_{\e,\pm}w_{\e,\pm}dx\\
&\rightarrow\int_{\R^{N+1}}(|\nabla U|^2+U^2)dx-\int_{\R^{N+1}}\mu f(U)Udx\\
&=(1-\mu)\int_{\R^{N+1}}f(U)Udx,
\endaligned
$$
which implies that there exist $C>0, t_\pm\in(0,1)$ and $\e_0>0$ such that
\begin{equation}\label{eqn:2order}
g''_{\e,\pm}(t)<-C\quad\text{for}\,\,\e\leq \e_0\,\,\text{and}\,\,t\in(1-t_\pm,1+t_\pm).
\end{equation}
By (\ref{eqn:zhengfu}), we have $g'_{\e,\pm}(t_{\e,\pm})=J'(t_{\e,\pm}w_{\e,\pm})w_{\e,\pm}=0$, and then the following holds for $\e$ small enough,
$$
J(w_{\e,\pm})=J(t_{\e,\pm}w_{\e,\pm})+\frac{1}{2}g''_{\e,\pm}(\xi_{\e,\pm})(t_{\e,\pm}-1)^2
$$
for $\xi_{\e,\pm}\in (1-t_\pm,1+t_\pm)$, which implies by (\ref{eqn:2order}), (\ref{eqn:decayv}), (\ref{eqn:decayV}) and the fact that
$v_{\e,\pm}\in
H_0^1(B_{\frac{R_{\e,\pm}}{\e}}(0))$ are the positive critical points of $J_{B_{\frac{R_{\e,\pm}}{\e}}(0)}$,
that for $\e$ small enough,
\begin{equation}\label{eqn:est1}
\aligned
&\frac{1}{2}C|t_{\e,\pm}-1|^2\\
\leq&J(t_{\e,\pm}w_{\e,\pm})-J(w_{\e,\pm})\\
=&\frac{|t_{\e,\pm}|^2-1}{2}\int_{|x|\leq\frac{R_{\e,\pm}}{\e}}(|\nabla v_{\e,\pm}|^2+|v_{\e,\pm}|^2)dx
-\int_{|x|\leq\frac{R_{\e,\pm}}{\e}}(F(t_{\e,\pm}v_{\e,\pm})-F(v_{\e,\pm}))dx\\
&-\frac{|t_{\e,\pm}|^2-1}{2}\int_{\mathcal{D}_1}(|\nabla v_{\e,\pm}|^2+|v_{\e,\pm}|^2)dx+\int_{\mathcal{D}_1}(F(v_{\e,\pm})-F(t_{\e,\pm}v_{\e,\pm}))dx\\
&+\frac{|t_{\e,\pm}|^2-1}{2}\int_{\mathcal{D}_2}(|\nabla V_{\e,\pm}|^2+|V_{\e,\pm}|^2)dx-\int_{\mathcal{D}_2}(F(V_{\e,\pm})-F(t_{\e,\pm}V_{\e,\pm}))dx\\
\leq& Ce^{-2(1-\delta)R_{\e,\pm}(\frac{1}{\e}-\delta')},
\endaligned
\end{equation}
where
$$
\aligned
&\mathcal{D}_1:=\bigg\{x\in\R^{N+1}|R_{\e,\pm}(\frac{1}{\e}-\delta')\leq|x|\leq\frac{R_{\e,\pm}}{\e}\bigg\},\\
&\mathcal{D}_2:=\bigg\{x\in\R^{N+1}||x|\geq R_{\e,\pm}(\frac{1}{\e}-\delta')\bigg\}.
\endaligned
$$
It follows from (\ref{eqn:est1}) that
\begin{equation}\label{eqn:est2}
|t_{\e,\pm}-1|\leq  Ce^{-(1-\delta)R_{\e,\pm}(\frac{1}{\e}-\delta')}.
\end{equation}
Let us note that
\begin{equation}\label{eqn:est3}
\aligned
J_{B_{\frac{R_{\e,\pm}}{\e}}(0)}(t_{\e,\pm}v_{\e,\pm})
= J(t_{\e,\pm}w_{\e,\pm}) +J_{\mathcal{D}_1}(t_{\e,\pm}v_{\e,\pm})-J_{\mathcal{D}_2}(t_{\e,\pm}V_{\e,\pm}),
\endaligned
\end{equation}
where $J_{\mathcal{D}_1}$, $J_{\mathcal{D}_2}$ denote
functional $J$ whose integral region are $\mathcal{D}_i, i=1,2$, respectively.
Combining (\ref{eqn:decayv}) and (\ref{eqn:decayV}), there exists $C>0$ such that for $t\in(\frac{1}{2},\frac{3}{2})$
$$
|J'_{\mathcal{D}_1}(tv_{\e,\pm})v_{\e,\pm}-J'_{\mathcal{D}_2}(tV_{\e,\pm})V_{\e,\pm}|\leq C e^{-2(1-\delta)R_{\e,\pm}(\frac{1}{\e}-\delta')},
$$
which implies by (\ref{eqn:est2}) that for $\e$ small enough
\begin{equation}\label{eqn:est4}
|J_{\mathcal{D}_1}(t_{\e,\pm}v_{\e,\pm})-J_{\mathcal{D}_2}(t_{\e,\pm}V_{\e,\pm})-J_{\mathcal{D}_1}(v_{\e,\pm})
+J_{\mathcal{D}_2}(V_{\e,\pm})|\leq C e^{-3(1-\delta)R_{\e,\pm}(\frac{1}{\e}-\delta')}.
\end{equation}
Using integration by parts, by ($f_1$), ($f_2$) and (\ref{eqn:decay0}) we have for any $\sigma\in(0,1)$
$$
\aligned
J_{\mathcal{D}_2}(V_{\e,\pm})=&\frac{2-\sigma}{2}J'_{\mathcal{D}_2}(V_{\e,\pm})V_{\e,\pm}
-\frac{1-\sigma}{2}\int_{\mathcal{D}_2}(|\nabla V_{\e,\pm}|^2+|V_{\e,\pm}|^2)dx\\
&+\int_{\mathcal{D}_2}(\frac{2-\sigma}{2}f(V_{\e,\pm})V_{\e,\pm}-F(V_{\e,\pm}))dx\\
\leq&-\frac{2-\sigma}{2}\omega_N[R_{\e,\pm}(\frac{1}{\e}-\delta')]^{N}V_{\e,\pm}(R_{\e,\pm}(\frac{1}{\e}-\delta'))
\frac{dV_{\e,\pm}(R_{\e,\pm}(\frac{1}{\e}-\delta'))}{dr},
\endaligned
$$
where we also use the fact that
$$
\int_{\mathcal{D}_2}[-V_{\e,\pm}\triangle V_{\e,\pm}+|V_{\e,\pm}|^2-f(V_{\e,\pm})V_{\e,\pm}]dx=0.
$$
Using the similar argument, we also deduce that for any $\sigma\in(0,1)$,
$$
J_{\mathcal{D}_1}(v_{\e,\pm})\geq\frac{-\sigma}{2}\omega_N[R_{\e,\pm}(\frac{1}{\e}-\delta')]^{N}
v_{\e,\pm}(R_{\e,\pm}(\frac{1}{\e}-\delta'))\frac{dv_{\e,\pm}(R_{\e,\pm}(\frac{1}{\e}-\delta'))}{dr}.
$$
According to the above facts, we have
\begin{equation}\label{eqn:est5}
\aligned
&J_{\mathcal{D}_1}(v_{\e,\pm})-J_{\mathcal{D}_2}(V_{\e,\pm})\\
\geq& \omega_N[R_{\e,\pm}(\frac{1}{\e}-\delta')]^{N}v_{\e,\pm}(R_{\e,\pm}(\frac{1}{\e}-\delta'))\times\\
&\bigg(\frac{2-\sigma}{2}\frac{d[V_{\e,\pm}(R_{\e,\pm}(\frac{1}{\e}-\delta'))]}{dr}
-\frac{\sigma}{2}\frac{d[v_{\e,\pm}(R_{\e,\pm}
(\frac{1}{\e}-\delta'))]}{dr}\bigg).
\endaligned
\end{equation}
An easy computation yields
$$
-\frac{d^2\tilde{V}_{\e,\pm}}{dr^2}+\frac{N}{r+R_{\e,\pm}
(\frac{1}{\e}-\delta')}\frac{d\tilde{V}_{\e,\pm}}{dr}+ \tilde{V}_{\e,\pm}=\frac{f(V_{\e,\pm}(R_{\e,\pm}
(\frac{1}{\e}-\delta'))\tilde{V}_{\e,\pm})}{V_{\e,\pm}(R_{\e,\pm}
(\frac{1}{\e}-\delta'))},\quad r>0
$$
and
$$
\left\{
  \begin{array}{ll}
  \displaystyle-\frac{d^2\tilde{v}_{\e,\pm}}{dr^2}+\frac{N}{r+R_{\e,\pm}
(\frac{1}{\e}-\delta')}\frac{d\tilde{v}_{\e,\pm}}{dr}+ \tilde{v}_{\e,\pm}=\frac{f(v_{\e,\pm}(R_{\e,\pm}
(\frac{1}{\e}-\delta'))\tilde{v}_{\e,\pm})}{v_{\e,\pm}(R_{\e,\pm}
(\frac{1}{\e}-\delta'))},&\quad r\in (0,R_{\e,\pm}\delta'),\\
 \tilde{v}_{\e,\pm}(R_{\e,\pm}\delta')=0,
  \end{array}
\right.
$$
when
$$
\tilde{V}_{\e,\pm}(x)=V_{\e,\pm}(|x|+R_{\e,\pm}(\frac{1}{\e}-\delta'))/V_{\e,\pm}(R_{\e,\pm}(\frac{1}{\e}-\delta'))
$$
and
$$
\tilde{v}_{\e,\pm}(x)=v_{\e,\pm}(|x|+R_{\e,\pm}(\frac{1}{\e}-\delta'))/v_{\e,\pm}(R_{\e,\pm}(\frac{1}{\e}-\delta')),
$$
respectively. Using ($f_1$) and (\ref{eqn:decay0}), and the standard elliptic estimates, we have that
$\tilde{V}_{\e,\pm}\rightarrow \tilde{V}_{\pm}$ in $C^2_{loc}[0,\infty)$ as $\e\rightarrow0^+$, and $\tilde{V}_{\pm}$ satisfies
$$
\frac{d^2\tilde{V}_{\pm}}{dr^2}-\tilde{V}_{\pm}=0,\quad r\in(0,\infty),\,\,\tilde{V}_{\pm}(0)=1,\quad \tilde{V}_{\pm}(r)\leq1,
$$
and then $\tilde{V}_{\pm}(r)=e^{-r}$.
Similarly, $\tilde{v}_{\e,\pm}\rightarrow \tilde{v}_{\pm}$ in $C^2_{loc}[0,1)$ as $\e\rightarrow0^+$, and $\tilde{v}_{\pm}$ satisfies
$$
\frac{d^2\tilde{v}_{\pm}}{dr^2}-\tilde{v}_{\pm}=0,\quad r\in(0,1),\,\,\tilde{v}_{\pm}(0)=1,\quad \tilde{v}_{\pm}(1)=0.
$$
It is easy to see from the above facts that
$$
\frac{d(\tilde{V}_{\pm}(r)-\tilde{v}_{\pm}(r))}{dr}\bigg|_{r=0}>0,
$$
which implies by the uniform convergence of $\tilde{V}_{\e,\pm}$ and $\tilde{v}_{\e,\pm}$ in $C^2_{loc}[0,1)$ that
for $\e>0$ small enough
$$
\frac{d(V_{\e,\pm}(r)-v_{\e,\pm}(r))}{dr}\bigg|_{r=R_{\e,\pm}(\frac{1}{\e}-\delta')}\geq C v_{\e,\pm}(R_{\e,\pm}(\frac{1}{\e}-\delta'))).
$$
Let $\sigma$ close sufficiently to $1$ in (\ref{eqn:est5}), then by (\ref{eqn:decayv}) we have
\begin{equation}\label{eqn:est6}
\aligned
&J_{\mathcal{D}_1}(v_{\e,\pm})-J_{\mathcal{D}_2}(V_{\e,\pm})\\
&\geq c(\pi)[R_{\e,\pm}(\frac{1}{\e}-\delta')]^{N}\cdot[v_{\e,\pm}(R_{\e,\pm}(\frac{1}{\e}-\delta'))]^2\\
&\geq Ce^{-2(1+\delta)R_{\e,\pm}(\frac{1}{\e}-\delta')}
\endaligned
\end{equation}
for some $C>0$.
It then follows from (\ref{eqn:est4}) that there exists $C>0$ such that
\begin{equation}\label{eqn:est7}
\aligned
J_{\mathcal{D}_1}(t_{\e,\pm}v_{\e,\pm})-J_{\mathcal{D}_2}(t_{\e,\pm}V_{\e,\pm})&\geq J_{\mathcal{D}_1}(v_{\e,\pm})
-J_{\mathcal{D}_2}(V_{\e,\pm})- C e^{-3(1-\delta)R_{\e,\pm}(\frac{1}{\e}-\delta')}\\
&\geq Ce^{-2(1+\delta)R_{\e,\pm}(\frac{1}{\e}-\delta')},
\endaligned
\end{equation}
which, together with (\ref{eqn:est3}), implies that the conclusion of Lemma \ref{Lem:estimate} holds.
The proof is complete.
\qed

\begin{lemma}\label{Lem:2lower}
For $\e>0$ small enough, we have
$$
c_\e\geq \e^{N+1}\bigg(4 J(U)+Ce^{-2(1+\delta)R_{\e,+}(\frac{1}{\e}-\delta')}+Ce^{-2(1+\delta)R_{\e,-}(\frac{1}{\e}-\delta')}\bigg)
$$
with $\delta>0$ small enough and independently of $\e$ and $\delta'$ has been given in (\ref{eqn:2-5-0}).
\end{lemma}
\Proof
Since $v_{\e,\pm}$  are positive ground state solutions of equation (\ref{eqn:half-eq1,2}), it is easy to see that
$$
J_{B_{\frac{R_{\e,\pm}}{\e}}(0)}(t_{\e,\pm}v_{\e,\pm})\leq J_{B_{\frac{R_{\e,\pm}}{\e}}(0)}(v_{\e,\pm}).
$$
Thus,
it follows from (\ref{eqn:zhengfu2}) that
\begin{equation}\label{eqn:zhengfu3}
\aligned
J_{B_{\frac{R_{\e,\pm}}{\e}}(0)}(v_{\e,\pm})
&\geq J(t_{\e,\pm}w_{\e,\pm}) +Ce^{-2(1+\delta)R_{\e,\pm}(\frac{1}{\e}-\delta')}\\
&\geq J(U) +Ce^{-2(1+\delta)R_{\e,\pm}(\frac{1}{\e}-\delta')}.
\endaligned
\end{equation}
On the other hand, we take $t_{\e,i},s_{\e,i}>0, i=1,2$ in (\ref{eqn:2-5}) such that
$$
J'_{\e,B_{\frac{R_{\e,+}}{\e}}(0)}(t_{\e,i}v_{i,\e}^+)v_{i,\e}^+=0,
\quad J'_{\e,B_{\frac{R_{\e,-}}{\e}}(0)}(s_{\e,i}v_{i,\e}^-)v_{i,\e}^-=0.
$$
Therefore, by the fact that $v_{\e,\pm}$ are positive radial ground state solutions of equation (\ref{eqn:half-eq1,2}),
and combining (\ref{eqn:2-5}) and (\ref{eqn:zhengfu3}), we have the conclusion of Lemma \ref{Lem:2lower}.
\qed

\emph{The proof of Theorem \ref{Thm:middle}.}
Let
$$
\psi(P,Q):=\min\{\frac{1}{2}|P-Q|,dist(P,\partial\Omega),dist(Q,\partial\Omega)\},
$$
and
$$\mathcal{F}_{0}=\max\limits_{(P,Q)\in\mathcal{L}^1\times\mathcal{L}^1}\psi(P,Q).
$$
It follows from Lemma \ref{Lem:2upper} and Lemma \ref{Lem:2lower} that
$$
(1-\delta)\mathcal{F}_{a_1,a_2}\leq (1+\delta)\lim\limits_{\e\rightarrow0^+}\min\{R_{\e,+},R_{\e,-}\}(1-\e\delta').
$$
Since $\delta'$ can be taken small arbitrarily and $a_1,a_2$ are arbitrary, we can obtain that
\begin{equation}\label{eqn:zhengfu4}
0<\mathcal{F}_{0}\leq\lim\limits_{\e\rightarrow0^+}\min\{R_{\e,+},R_{\e,-}\}.
\end{equation}
On the other hand, since maximum and minimum points $P_\e,Q_\e$ belong to the $x_{N+1}$-axis for sufficiently small $\e$,
there exist point $\mathcal{R}_\e\in \overline{P_\e Q_\e}$ such that $\mathcal{R}_\e\in\partial\Omega_{\e,\pm}^1$. Moreover,
 for some $b_{\e,\pm}>0$,
$$
|P_\e-\mathcal{R}_\e|=\frac{b_{\e,+}}{b_{\e,+}+b_{\e,-}}|P_\e-Q_\e|,\quad|Q_\e-\mathcal{R}_\e|=\frac{b_{\e,-}}{b_{\e,+}+b_{\e,-}}|P_\e-Q_\e|.
$$
Thus, based on the definition of $R_{\e,\pm}$, we have immediately
$$
R_{\e,\pm}\leq\min\bigg\{\frac{b_{\e,+}}{b_{\e,+}+b_{\e,-}}|P_\e-Q_\e|,\frac{b_{\e,-}}{b_{\e,+}+b_{\e,-}}|P_\e-Q_\e|,
dist(P_\e,\partial\Omega),dist(Q_\e,\partial\Omega)\bigg\}
$$
By (\ref{eqn:zhengfu4}) and the definition of $\mathcal{F}_{0}$, we obtain
$$
\frac{b_{\e,\pm}}{b_{\e,+}+b_{\e,-}}\rightarrow\frac{1}{2}
$$
as $\e\rightarrow0^+$. That is to say, $\psi(P_\e,Q_\e)\rightarrow\mathcal{F}_0$ as $\e\rightarrow0^+$,
and by (\ref{eqn:limit}), one has
$\psi(\bar{P},\bar{Q})=\mathcal{F}_0$. Since nodal solution $u_\e$ of equation (\ref{eqn:redu-qs-})
belonges to $ H_{\sharp}(\Omega)$, we also have
$\psi(-P_\e,-Q_\e)\rightarrow\mathcal{F}_0$ as $\e\rightarrow0^+$, and
$\psi(-\bar{P},-\bar{Q})=\mathcal{F}_0$.
Since maximum and minimum points of $u_\e$ belong to the $x_{N+1}$-axis for sufficiently small $\e$,
by using the similar arguments as Lemma 5.1 in \cite{Noussair97}, we can obtain
$$
\aligned
&dist(\bar{P},\partial\Omega^1)=dist(\bar{Q},\partial\Omega^1)=\frac{1}{2}|\bar{P}-\bar{Q}|,\\
&dist(-\bar{P},\partial\Omega^2)=dist(-\bar{Q},\partial\Omega^2)=\frac{1}{2}|\bar{P}-\bar{Q}|,
\endaligned
$$
where $\Omega^1,\Omega^2$ have been defined in (\ref{eqn:omega+}). It is easy to conclude from the above identity that
the locations of four concentration points are as follows
$$
\bar{P}=(\mathbf{0},\frac{3a^2+b^2}{8}),\quad \bar{Q}=(\mathbf{0},\frac{a^2+3b^2}{8}),
\quad-\bar{P}=(\mathbf{0},-\frac{3a^2+b^2}{8}),\quad -\bar{Q}=(\mathbf{0},-\frac{a^2+3b^2}{8}),
$$
or,
$$
\bar{P}=(\mathbf{0},\frac{a^2+3b^2}{8}),\quad \bar{Q}=(\mathbf{0},\frac{3a^2+b^2}{8}),
\quad-\bar{P}=(\mathbf{0},-\frac{a^2+3b^2}{8}),\quad -\bar{Q}=(\mathbf{0},-\frac{3a^2+b^2}{8}).
$$
Recalling Lemma \ref{Lem:mox} and Theorem \ref{Thm:exist}, we have proved that there exists a nonradial
nodal solution $u_\e$ of equation (\ref{eqn:redu-qs-}) whose maximum and minimum points concentrate
 at four points in the $x_{N+1}$-axis contained in $\Omega$.
 The corresponding solution of equation (\ref{eqn:qs}), still denoted by $u_\e$, concentrates exactly at four
 $(N-1)$-dimensional spheres in $A$, two of them (corresponding to $\bar{P}$, $\bar{Q}$)
 are characterized by the angle $\theta=0$, and by the radial coordinate
 $$
 r=\frac{1}{2}\sqrt{3a^2+b^2},\quad\text{and}\quad r=\frac{1}{2}\sqrt{a^2+3b^2}
 $$
 the others (corresponding to $-\bar{P}$, $-\bar{Q}$) paced the angle $\theta=\pi$ and the radial coordinate
 $$
 r=\frac{1}{2}\sqrt{a^2+3b^2},\quad\text{and}\quad r=\frac{1}{2}\sqrt{3a^2+b^2}.
 $$
The proof is complete.
\qed

\end{document}